\documentclass[a4paper,fleqn]{cas-sc}

\usepackage[numbers]{natbib}
\usepackage{amssymb}
\usepackage{amsmath}
\usepackage{siunitx}
\usepackage[algoruled, linesnumbered, lined]{algorithm2e} % load before cleveref
\usepackage[capitalise]{cleveref}
\usepackage{booktabs}
\usepackage{subcaption}
\usepackage{mathtools} % \coloneqq 
\usepackage{bm}

\usepackage{xcolor}
\definecolorset{RGB}{tol}{}{blue,68,119,170;cyan,102,204,238;green,34,136,51;yellow,204,187,68;red,238,102,119;purple,170,51,119;grey,187,187,187}

\newif\ifbuild
\buildfalse

\ifbuild
\usepackage{tikz}
\usepackage{pgfplots}
\usepackage{pgfplotstable}

\usetikzlibrary{positioning, math,calc, backgrounds, fit}
\usepgfplotslibrary{statistics,groupplots}

\usetikzlibrary{external}
\pgfplotscreateplotcyclelist{bright}{%
	tolblue,every mark/.append style={fill=tolblue!80!black},mark=*\\%
	tolcyan,every mark/.append style={fill=tolcyan!80!black},mark=square*\\%
	tolgreen,every mark/.append style={fill=tolgreen!80!black},mark=triangle*\\%
	tolyellow,every mark/.append style={fill=tolyellow!80!black},mark=star\\%
	tolred,every mark/.append style={fill=tolred!80!black},mark=diamond*\\%
	tolpurple,every mark/.append style={solid,fill=tolpurple!80!black},mark=pentagon*\\%
	tolgrey,every mark/.append style={solid,fill=tolgrey!80!black},mark=10-pointed star\\%
}

\pgfplotscreateplotcyclelist{bright-grouped}{%
	tolblue,every mark/.append style={fill=tolblue!80!black},mark=*\\% coarse
	tolcyan,every mark/.append style={fill=tolcyan!80!black},mark=square*\\% fine
	tolgreen,every mark/.append style={fill=tolgreen!80!black},mark=triangle*\\%
	tolgreen,every mark/.append style={fill=tolgreen!80!black},mark=triangle*\\%
	tolgreen,every mark/.append style={fill=tolgreen!80!black},mark=triangle*\\%
	tolgreen,every mark/.append style={fill=tolgreen!80!black},mark=triangle*\\%
	tolyellow,every mark/.append style={fill=tolyellow!80!black},mark=star\\%
	tolyellow,every mark/.append style={fill=tolyellow!80!black},mark=star\\%
	tolyellow,every mark/.append style={fill=tolyellow!80!black},mark=star\\%
	tolyellow,every mark/.append style={fill=tolyellow!80!black},mark=star\\%
	tolred,every mark/.append style={fill=tolred!80!black},mark=diamond*\\%
	tolred,every mark/.append style={fill=tolred!80!black},mark=diamond*\\%
	tolred,every mark/.append style={fill=tolred!80!black},mark=diamond*\\%
	tolred,every mark/.append style={fill=tolred!80!black},mark=diamond*\\%
	tolpurple,every mark/.append style={fill=tolpurple!80!black},mark=pentagon*\\%
	tolpurple,every mark/.append style={fill=tolpurple!80!black},mark=pentagon*\\%
	tolpurple,every mark/.append style={fill=tolpurple!80!black},mark=pentagon*\\%
	tolpurple,every mark/.append style={fill=tolpurple!80!black},mark=pentagon*\\%
}

\pgfplotsset{compat=1.18,
	/pgfplots/bar  cycle  list/.style={/pgfplots/cycle  list={%
			{tolblue!20!black,fill=tolblue,mark=none},%
			{tolcyan!20!black,fill=tolcyan,mark=none},%
			{tolgreen!20!black,fill=tolgreen,mark=none},%
			{tolyellow!20!black,fill=tolyellow,mark=none},%
			{tolred!20!black,fill=tolred,mark=none},%
			{tolpurple!20!black,fill=tolpurple,mark=none},%
			{tolgrey!20!black,fill=tolgrey,mark=none},%
		}
	},
	cycle list name=bright,
	only if/.style args={entry of #1 is #2}{
		/pgfplots/boxplot/data filter/.code={
			\edef\tempa{\thisrow{#1}}
			\edef\tempb{#2}
			\ifx\tempa\tempb
			\else
			
			\fi
		}
	},
	table/search path={data},
    layers/my layer set/.define layer set={
            background,
            main,
            foreground
    }{
    },
    set layers=my layer set,
}

\def\mystrut{{\vphantom{hg}}}

\pgfplotsset{
	legend image with text/.style={
		legend image code/.code={%
			\node[anchor=center] at (0.3cm,0cm) {#1};
		}
	},
}

\fi % end ifbuild

\newcommand{\R}{\mathbb{R}}
\newcommand{\N}{\mathbb{N}}
\newcommand{\refsol}[1]{#1}
\newcommand{\nnsol}[1]{\tilde{#1}}
\newcommand{\intsol}[1]{\hat{#1}}
\newcommand{\Rey}{\Re} % \mathrm{Re}
\newcommand{\differential}{\,\text{d}}
\renewcommand{\vec}[1]{\mathbf{#1}}
\newcommand{\fvec}[1]{\bm{#1}}
\newcommand{\ts}{\Delta t}

\usepackage{lineno}

\begin{document}
\let\WriteBookmarks\relax
\def\floatpagepagefraction{1}
\def\textpagefraction{.001}

\shorttitle{}

\shortauthors{R. Jendersie et~al.} 

%\begin{frontmatter}



\title [mode = title] {A robust and stable hybrid neural network/finite element method for 2D flows that generalizes to different geometries}

\author[ian,isg]{Robert Jendersie}\ead{robert.jendersie@ovgu.de} % [orcid=0000-0001-8693-0222]
\cormark[1]
\author[ian]{Nils Margenberg}\ead{nils.margenberg@ovgu.de} % [orcid=0000-0003-2089-5545]
\author[ecmwf,isg]{Christian Lessig}\ead{christian.lessig@ovgu.de} % [orcid=0000-0002-2740-6815]
\author[ian]{Thomas Richter}\ead{thomas.richter@ovgu.de} % [orcid=0000-0003-0206-3606]

\cortext[cor1]{Corresponding author}
\affiliation[ian]{organization={Institute of Analysis and Numerics, Otto-von-Guericke Universität},%Department and Organization
            addressline={Universitätsplatz 2}, 
            city={Magdeburg},
            postcode={39104}, 
            country={Germany}}

\affiliation[isg]{organization={Institute of Simulation and Graphics, Otto-von-Guericke Universität},
            addressline={Universitätsplatz 2}, 
            city={Magdeburg},
            postcode={39104}, 
            country={Germany}}
\affiliation[ecmwf]{organization={European Centre for Medium-Range Weather Forecasts},
            addressline={Robert-Schuman-Platz 3}, 
            city={Bonn},
            postcode={53175}, 
            country={Germany}}

%% Abstract
\begin{abstract}
The deep neural network multigrid solver (DNN-MG) combines a coarse-grid finite element simulation with a deep neural network that corrects the solution on finer grid levels, thereby improving the computational efficiency. 
%By only making predictions on small patches of the mesh, the neural network learns a local correction which allows it to generalize beyond the training domains. 
In this work, we discuss various design choices for the DNN-MG method and demonstrate significant improvements in accuracy and generalizability when applied to the solution of the nonstationary Navier-Stokes equations.
We investigate the stability of the hybrid simulation and show how the neural networks can be made more robust with the help of replay buffers. 
After an initial single-step training, we run the hybrid simulation for extended periods and compute new reference solutions of the neural network perturbed state for each step. By retraining on this data, the error caused by the neural network over multiple time-steps due to distributional shift can be effectively reduced without the need for a differentiable numerical solver.
Furthermore, we compare multiple neural network architectures, including recurrent neural networks and Transformers, and study their ability to utilize more information from an increased temporal and spatial receptive field. Transformers allow us to make use of information from cells outside the predicted patch even with unstructured meshes while maintaining the locality of our approach. This can further improve the accuracy of DNN-MG without a significant impact on performance.
\end{abstract}



%% Keywords
\begin{keywords}
Computational fluid dynamics \sep 
Navier-Stokes equations \sep 
Finite element method \sep
Deep learning \sep
Stability


\end{keywords}

\maketitle

%\end{frontmatter}

%%%%%%%%%%%%%%%%%%%%%%%%%%%%%%%%%%%%%%%%%%%%%%%%%%%%%%%%%%%%%%%%%%%%%%%%%%%%
\section{Introduction}
The accurate solution of the nonstationary Navier-Stokes equations is of great importance in engineering and science. 
Existing numerical methods can provide accurate and reliable solutions. 
However, at high resolutions or Reynolds numbers, these methods become prohibitively expensive. 

Finite element methods are among the most successful approaches for the solution of partial differential equations. For fluid flow problems, existing methods such as multigrid methods are highly efficient and can achieve optimal computational complexity~\cite{kimmritzParallelMultigrid2011}. 
The majority of the costs is thereby typically in the linear solver.
Consequently, choosing the right solver for a specific problem is vital for computational efficiency~\cite{ahmedAssessmentSolversSaddle2018, Ghai2018}. Further work to speed up such numerical methods focuses on the use of massively parallel hardware such as graphics processing units (GPUs). However, sparse solvers generally achieve only a low utilization of GPUs due to memory bandwidth limitations~\cite{Liebchen2024, Thomas2024}. Currently, matrix free finite element approaches promise the best potential for high computational efficiency on massively parallel computers~\cite{Munch2023}. However, the challenges of GPU implementations are considerable, especially if a high degree of flexibility is to be maintained with regard to the numerical methodology and the type of PDE taken into account. 

In recent years, data-driven methods using machine learning have emerged as an alternative to conventional approaches. Deep neural networks (NNs) open new avenues for the solution of PDEs and, having co-evolved with GPUs, can effectively leverage the massive compute that they provide. Neural networks can directly learn a specific solution by using the differential equation as loss function, completely sidestepping the need for a discretization. First introduced in~\cite{Lagaris1998}, this idea has garnered increased interest since its reemergence as Physics-informed neural networks (PINNs)~\cite{Raissi2019a, eDeepRitzMethod2018, luDeepXDEDeepLearning2021}. The approach is well suited for the solution of inverse problems but is usually not competitive with traditional numerical methods for forward problems, which is the setting we are interested in~\cite{Grossmann2024,Toscano2025}. See~\cite{Toscano2025} for a recent survey.

The efficiency of neural network-based methods can be substantially improved by learning a solution operator instead of a single solution. Finite-dimensional operators can be learned from fully discretized inputs and outputs~\cite{Guo2016, Bhatnagar2019,eichingerStationaryFlowPredictions2021,eichinger2022}. DeepONets~\cite{luLearningNonlinearOperators2021} do away with the need for a discretization of the outputs and Neural Operators~\cite{JMLR:v24:21-1524} are also flexible with regard to the representation of the inputs. 
A challenge for learned operators is to extrapolate beyond the training interval for time-dependent problems~\cite{Nayak2025}. Recent work improves upon this by introducing an auto-regressive pathway that incorporates operator evaluations from past time-points~\cite{Diab2025} or by learning to predict the time-derivative and incorporating the operator with a Runge-Kutta scheme~\cite{Nayak2025}.
 
To get accurate solutions over longer time intervals, time-stepping remains the standard method even with neural networks. Together with a full discretization in space and time, neural networks can be used in place of a conventional scheme or replace parts of it. By only changing heuristics, such as the preconditioner of a linear solver with machine learning, guaranties on convergence and accuracy can be preserved. However, the potential for speedup is also limited~\cite{luzLearningAlgebraicMultigrid2020, Markidis2021, huangLearningOptimalMultigrid2023, Holguin2024}. 

Greater gains are possible with super-resolution methods, where a numerical method may be used to cheaply compute a low-fidelity solution, which is then enhanced by a deep neural network; see \cite{Sofos2025} for a recent survey of this approach for fluid flows. In computer vision and graphics, this process is commonly called upscaling, while the weather and climate community refers to high-resolution reconstruction as downscaling. Generally these terms refer to just one direction of information flow, from the numerical method to the neural network. The task is therefore to reconstruct a statistically consistent, temporally coherent high-resolution realization~\cite{Xie2018, Schmidt2025}. A tighter integration with the numerical method is needed when the neural network has to correct dynamics that only develop with higher resolutions and which can effect the coarse-scale flow. This is accomplished either by directly updating the coarse solution, possibly without explicitly forming a fine solution, during the time-stepping~\cite{Kochkov2024, Gregory2025, Witte2025}, or by introducing the neural network corrections through a forcing term~\cite{fabraFiniteElementApproximation2022,delaraAcceleratingHighOrder2022,manriquedelaraAcceleratingHighOrder2023}.

Another approach to time-stepping with neural networks keeps only the discretization and learns the full update from data. In~\cite{Brandstetter2022} the authors present a general scheme that mimics the structure of conventional stencil based schemes but does every computation with a graph neural network. In weather forecasting, fully data-driven models are already operational~\cite{biPanguWeather3DHighResolution2022,Lam2023,Lang2024} as there are wealth of historical data to learn from and important physical processes that cannot be easily modeled or resolved by numerical models. Nevertheless, NeuralGCM~\cite{Kochkov2024} demonstrates that a conventional solver enhanced by neural network corrections can be competitive.

Hybrid and fully data-driven methods that use a neural network auto-regressively , i.e. in a time-stepping loop, face similar challenges as conventional methods regarding stability. Over multiple time-steps, the error introduced by the neural network can cause a distributional shift in the network input data. This discrepancy between ideal training samples and network input at later time-steps during inference can lead to unphysical predictions and instabilities. To alleviate this, methods aim to directly minimize the error across many time-steps in addition to the local single-step error through unrolled training. In a hybrid setting, a fully differentiable code makes it possible to compute gradients of short simulation runs via backpropagation through time (BPTT)~\cite{Kochkov2024}. Even without BPTT, a multi-step roll-out can be helpful for stabilization in combination with a multi-step loss function~\cite{Brenowitz2018, Brandstetter2022}. Both approaches are compared in~\cite{List2025}. A simpler approach is to store single- or multi-step step training predictions in a replay buffer, to serve as new inputs in further training steps when combined with existing reference outputs~\cite{Chen2025, Chen2025a}.
Alternatively, completely new training samples can be generated by running a reference simulation in concert with the neural network to get corrected outputs for the input state perturbed by the neural network~\cite{Rasp2020}.
No modification to the training is necessary if the process is done by alternating data generation and training stages, but multiple iterations can be necessary to achieve the best results~\cite{Gregory2024, Gregory2025}. A completely different approach, named thermalization, is taken in \cite{Pedersen2025}, where a separate diffusion network is trained on the target distribution to stabilize the predictions of the primary neural network during inference.

The deep neural network multigrid solver (DNN-MG) is a hybrid time-stepping method, where the neural network explicitly reconstructs a fine solution, which uses corrective forcing in the high-fidelity space as feedback for the numerical simulation. The method was introduced in \cite{MargenbergRichter2020}, demonstrating its application to the 2D Navier-Stokes equations. On the same problem, we also investigated the use of the stream function to ensure divergence-freedom by construction~\cite{margenbergStructurePreservationDeep2021}. In more recent work, we showed the effectiveness of DNN-MG for the solution of the 3D Navier-Stokes equations~\cite{Margenberg2024}. The tight integration of a neural network with the finite element method has also enabled first steps towards an error analysis of our method~\cite{kapustsinHybridFiniteElement,kapustsinErrorAnalysisHybrid2023}.

In this work, we present improvements to the DNN-MG method that substantially strengthen its practical applicability.
\begin{figure}
\newlength{\channelHeight}
\setlength{\channelHeight}{2.15cm}
	\begin{subfigure}{0.39\textwidth}
    \centering
    \ifbuild
	\tikzsetnextfilename{ro1_comparison}
	\begin{tikzpicture}[text=black]
		\node(Img){\includegraphics[width=0.8\textwidth, trim={0 0 32cm 0}
			, clip]{instability.png}};
		\node(CoarseLabel)[anchor=north east, inner sep=0.25cm] at (Img.north east){MG(5)};
		\node(NNLabel)[anchor=north east, inner sep=0.25cm]  at ($(Img.north east) - (0,\channelHeight)$) {DNN-MG(5+1)};
		\node(FineLabel)[anchor=north east, inner sep=0.25cm] at ($(Img.north east) - (0, 2.0 * \channelHeight)$){MG(6)};
	\end{tikzpicture}
    \else
	\includegraphics{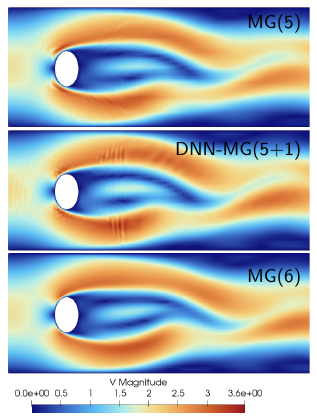}
	\fi
	\subcaption{Velocity magnitude at $t=\qty{1.17}{s}$}
	\label{fig:solution-nn-artifacts}
	\end{subfigure}
	\begin{subfigure}{0.59\textwidth}
        \ifbuild
    	\tikzsetnextfilename{ro1_divergence}
		\begin{tikzpicture}
			\begin{axis}[
				xlabel={time $t [\unit{s}]$},
				ylabel={$ J_{\text{div}}$},
				width=0.99\textwidth,
                height=7.4cm,
				ymin=0.0,
				ymax=6.2,
				]
				\addplot+[mark repeat=32] table[x index=0, y index=1, col sep=space]{fcts_gc1_coarse.txt};
				\addplot+[mark repeat=32] table[x index=0, y index=1, col sep=space]{fcts_gc1_ifix_instability_1_ref_detailed.txt};
				\addplot+[mark repeat=32] table[x index=0, y index=1, col sep=space]{fcts_gc1_fine.txt};
				\addplot +[mark=none,color=gray] coordinates {(1.17, 0) (1.17, 6.2)};
				\legend{MG(5), DNN-MG(5+1), MG(6)};
			\end{axis}
		\end{tikzpicture}
        \else
        \includegraphics{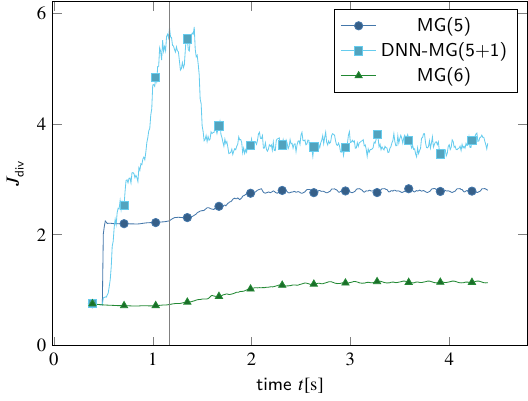}
        \fi
		\subcaption{Divergence over time}
		\label{fig:div-nn-artifacts}
	\end{subfigure}
	\caption{Comparison of the single obstacle case with reference simulations refinement levels 5 and 6 and the NN enhanced simulation using a small MLP to correct the solution on level 5. In all cases, the solution at level 6 is used until $t=0.5$. The results show that the DNN-MG method can develop instabilities that manifest as nonphysical high frequency features in the flow and that, by some metrics, make the solution worse than the level 5 baseline.}
    \label{fig:motivation_instability}
\end{figure}
One particular motivation is shown in~\cref{fig:motivation_instability} where the classical channel flow around an obstacle is considered. 
%Here, the DNN-MG method is trained on a different obstacle shape than used at inference time. 
The visualization on the left reveals spurious patterns in the DNN-MG solution and the divergence functional on the right shows their rapid formation from the accumulation of systemic errors.  
In the present work, we show that this instability is not inherent to the DNN-MG method but arises from the distributional shift between training and inference.
With this insight, we avoid the instability using replay buffers and training in multiple stages. After the initial training on samples derived from the numerical solver, additional data is generated with the hybrid simulation to condition the neural networks on the errors that they themselves introduce.

In our previous works, replay training was not needed. These stability issues emerge from further tuning of the neural network component of DNN-MG and replay training helps to translate improved accuracy on the training task into more accurate hybrid simulations.
%As replay training helps bridge the gap between the training task and inference, we achieve further improvements in simulation accuracy by tuning the neural network component of DNN-MG.
%Stabilizing training is necessary to translate the improved training loss from modifications described in this work into increased accuracy during inference. 
In particular, we consider different neural network designs and their ability to leverage additional spatial and temporal context to generalize to different geometries. We find that enlarging the spatial receptive field of the neural network is effective to improve the skill of the simulations without negatively affecting speed.

The remainder of the paper is structured as follows.
In \cref{sec:numerics} we describe the discretization and numerical solution of the Navier-Stokes equations, followed by its modifications for DNN-MG in~\cref{sec:dnn-mg}. \cref{sec:cases} gives details of the test cases we use and introduces suitable metrics for the evaluation of DNN-MG. Stabilization methods are discussed in \cref{sec:stability} and the generalization ability of DNN-MG is demonstrated in \cref{sec:generalization}. Finally, we summarize our findings in \cref{sec:conclusion}.

%%%%%%%%%%%%%%%%%%%%%%%%%%%%%%%%%%%%%%%%%%%%%%%%%%%%%%%%%%%%%%%%%%%%%%%%%%%%
\section{Numerical method}\label{sec:numerics}
The deep neural network multigrid solver (DNN-MG) is a hybrid numerical scheme for the solution of partial differential equations that combines a numerical solver with a deep neural network. We build on top of the finite element library Gascoigne~\cite{Gascoigne3D}, which uses the Newton-Krylov method with a geometric multigrid preconditioner with a cell-based Vanka-smoother. However, another choice for the numerical solver is possible. What is needed from the geometric multigrid toolkit for DNN-MG is a hierarchical mesh with compatible finite element spaces defined on at least two levels.

In this work, we consider the incompressible, nonstationary Navier-Stokes equations
\begin{equation}
  \begin{alignedat}{2}\label{eq:nsstrong}
    \partial_t \vec v + (\vec v\cdot \nabla)\vec v - \nu\Delta \vec v
    +\nabla p &= f \quad &&\text{on } [0,\,T] \times \Omega,\\
    \nabla \cdot \vec v &= 0 \quad &&\text{on } [0,\,T] \times \Omega,
  \end{alignedat}
\end{equation}
on a domain $\Omega \subset \R^2$ and the time interval $[0, T]$, where $\nu>0$ is the kinematic viscosity and $f:[0,T]\times\R^2 \to \R^2$ is the external forcing. We seek solutions for the velocity \(\vec v\colon [0,\,T]\times \Omega \to \R^2\) and pressure \(p\colon [0,\,T]\times \Omega \to \R\) subject to boundary conditions
\begin{equation}\label{eq:boundary}
  \begin{alignedat}{2}
    \vec v(0,\,\cdot) &= \vec v_0(\cdot)\quad &&\text{on }\Omega,\\
    \vec v &= \vec v^D \quad &&\text{on } [0,\,T] \times \Gamma^D,\\
    \nu(\vec n\cdot\nabla)\vec v - p\vec n &=0 \quad &&\text{in } [0,\,T] \times \Gamma^N,
  \end{alignedat}
\end{equation}
where \(\vec n\) denotes the outward facing unit normal vector of the boundary \(\partial\Omega\) of the domain. On the boundary \(\Gamma = \Gamma^D \cup \Gamma^N\) we prescribe Dirichlet boundary conditions \(\Gamma^D\) or Neumann boundary conditions \(\Gamma^N\).

%%%%%%%%%%%%%%%%%%%%%%%%%%%%%%%%%%%%%%%%%%%%%%%%%%%%%%%%%%%%%%%%%
\subsection{Variational formulation}
To apply the finite element method, we require the variational formulation of the Navier-Stokes equations.
Let \(L^2(\Omega)\) denote the space of square integrable functions on the
domain \(\Omega\subset\mathbb{R}^2\) with scalar product \((\cdot \, ,\cdot)\) and
let \(H^1(\Omega)\) be the subspace of \(L^2(\Omega)\) of functions with weak first derivatives in
\(L^2(\Omega)\). The function spaces for $\vec v$ and $p$ are then defined as
\begin{equation}\label{eq:VQ}
  \begin{alignedat}{1}
     \vec V &\coloneqq \vec v^D + H^1_0(\Omega;\Gamma^D)^2,\quad H_0^1 (\Omega;\Gamma^D)^2\coloneqq \left\{\vec v\in H^1(\Omega)^2\colon \vec v=0 \text{  on } \Gamma^D \right\}\\
    L &\coloneqq \left\{p\in L^2(\Omega),\text{ and, if }\Gamma^N=\emptyset,\; \int_{\Omega}p \differential x = 0\right\},
  \end{alignedat}
\end{equation}
where \(\vec v^D\in H^1(\Omega)^2\) is an extension of the Dirichlet data on \(\Gamma^D\) into the domain. If we have only Dirichlet boundaries, we normalize the pressure to yield uniqueness. With \eqref{eq:VQ}, the weak formulation of \eqref{eq:nsstrong} is given by
\begin{equation}
\begin{alignedat}{2}\label{eq:ns}
  (\partial_t \vec v,\,\fvec \phi) + (\vec v\cdot \nabla \vec v,\, \fvec \phi) +
  \nu(\nabla \vec v,\, \nabla \fvec\phi) -
  (p,\,\nabla\cdot \fvec\phi)
  &= (\vec f,\,\fvec \phi)\;\;&&\forall \fvec \phi \in H^1_0(\Omega;\Gamma^D)^2,\\
  (\nabla \cdot \vec v, \, \xi) &= 0\quad &&\forall \xi \in L,\\
  \vec v(0,\,\cdot ) &= \vec v_0(\cdot)\quad &&\text{on }\Omega.
\end{alignedat}
\end{equation}

%%%%%%%%%%%%%%%%%%%%%%%%%%%%%%%%%%%%%%%%%%%%%%%%%%%%%%%%%%%%%%%%%
\subsection{Discretization}\label{sec:discretization}
We discretize $\Omega$ using a quadrilateral finite element mesh $\Omega_h$ that satisfies the usual requirements on the structural and form regularity so that the standard interpolation results hold, cf.~\cite[Section 4.2]{richterFluidstructureInteractionsModels2017}. 
We use a standard $H^1$-conforming finite element space $Q^{(2)}_h$ of continuous isoparametric biquadratic elements on each quadrilateral $K \in \Omega_h$.
The trial- and test-functions are then defined as $\vec v_h,\,\fvec\psi_h \in \vec V_h = [Q^{(2)}_h]^2$ and $p_h, \xi_h \in L_h = Q^{(2)}_h$.
Since the resulting equal-order finite element pair \(\vec V_h\times L_h\) does not satisfy the inf-sup condition, we add stabilization terms based on local projections~\cite{beckerFiniteElementPressure2001}. The advantage of equal-order discretization lies in the simplified memory layout which allows to cluster all local degrees of freedom (for velocity and pressure) in the nodes~\cite{kimmritzParallelMultigrid2011}. Through interpolation of $\vec f$  and $\vec v^D$, we get the approximations $\vec f_h,\:\vec v_h^{D}\in C(I;\,\vec V_{h})$ to state the resulting semi-discrete variational problem as finding $\vec v_h$, $p_h$ such that $\vec v_{h}=\vec v_h^D\:\text{on}\:\bar{I}\times\Gamma_D$ and
\begin{equation}
  \begin{alignedat}{2}\label{eq:disc:ns}
    (\partial_t \vec v_h,\,\fvec\psi_h) + (\vec v_h\cdot \nabla \vec v_h,\, \fvec\psi_h)\qquad\\ +
    \nu (\nabla \vec v_h,\, \nabla \fvec\psi_h) -
    (p_h,\,\nabla\cdot \fvec\psi_h)
    &= (\vec f_h,\,\fvec\psi_h)\quad&\forall \fvec\psi_h\in  \vec V_h,\\
    (\nabla \cdot \vec v_h, \, \xi_h)
    +\sum_{K\in\Omega_h}\alpha_K (\nabla (p_h-\pi_h p_h),\nabla (\xi_h-\pi_h \xi_h))
    &= 0\quad &\forall \xi_h \in L_h\,.
  \end{alignedat}
\end{equation}
Here, $\pi_h:Q^{(2)}_h\to Q^{(1)}_h$ denotes the interpolation into the space $Q^{(1)}_h$ of isoparametric bilinear elements on each quadrilateral. The local stabilization parameter $\alpha_K$ is computed as
\[
 \alpha_K = \alpha_0 \left(\frac{\nu}{h_K^2}+\frac{\|\vec v_h\|_{\infty,K}}{h_K}\right)^{-1},
\]
where $h_K$ is diameter of element $K$ and $\alpha_0 > 0$ is a constant, set to $\alpha_0=0.04$ in all our numerical experiments.

In time, we discretize the interval $[0, T]$ into $N \in \N$ steps of uniform size
\[
0=t_0<t_1<\cdots <t_N=T,\: \ts = t_n-t_{n-1},
\]
leading to the fully discrete approximation of the velocity and pressure 
\[
\vec{v}_n \coloneqq \vec v_h(t_n), \qquad p_n\coloneqq p_h(t_n),
\]
at time $t_n$. We apply the second-order Crank-Nicolson scheme~\cite{SonnerRichter2017} to~\eqref{eq:disc:ns}, resulting in the fully discrete problem
\begin{equation}
\begin{alignedat}{4}\label{eq:4cranknicholson3}
\frac{1}{\ts}(\vec v_n,\,\fvec\phi_h)\,
+{}\frac{1}{2} (\vec v_n\cdot \nabla \vec v_n,\,\fvec\phi_h)
+\frac{\nu}{2}(\nabla \vec v_n,\,\nabla \fvec\phi_h)
-(p_n,\,\nabla \cdot \fvec\phi_h)
    \\ = \frac{1}{\Delta t}(\vec v_{n-1},\, \fvec\phi_h) +\frac{1}{2}(\vec f_n,\fvec\phi_h)
  +\frac{1}{2}(\vec f_{n-1},\,\fvec\phi_h) & \\
  -\frac{1}{2}(\vec v_{n-1}\cdot\nabla \vec v_{n-1},\,\fvec\phi_h)
  - \frac{\nu}{2}(\nabla \vec v_{n-1},\,\nabla \fvec\phi_h)
   & \quad \forall \fvec\phi_h\in \vec V_h, \\
%%%
(\nabla \cdot \vec v_{n}, \, \xi_h)+
\sum_{K\in\Omega_h}\alpha_K (\nabla (p_n-\pi_h p_n),\nabla (\xi_h-\pi_h \xi_h))
=0 & \quad \forall \xi_h \in L_h.
\end{alignedat}
\end{equation}
With the combined solution vector $\vec{x}_n = (p_n, \vec{v}_n)$ we shorten \cref{eq:4cranknicholson3} to
\begin{equation}\label{eq:nonlinearshort}
    \mathcal{A}_h(\vec{x}_n) = \mathcal{F}_h(\vec{v}_{n-1}, \vec{f}_{n-1}, \vec{f}_n),
\end{equation}
where $[\mathcal{A}_h(\vec x_n)]_i$ and $[\mathcal{F}_h]_i$ are, respectively, the left and right hand sides for all test functions $(\xi_h^i,\,\fvec \phi_h^i)$.

\subsection{Numerical Solver}
We solve the system of algebraic equations~\eqref{eq:nonlinearshort} with Newton's method. We ignore the velocity dependence of the stabilization term in \cref{eq:4cranknicholson3} when assembling the Jacobian, as is commonly done. The Jacobian matrix is only recomputed when the convergence rate
\[
\rho= \frac{\lVert{\vec r^{(j+1)}_n}\rVert}{\lVert{\vec r^{(j)}_n}\rVert}, \quad \vec{r}_n^{(i)} = \mathcal{F}_h(\vec{x}_{n-1}, \vec{f}_n, \vec{f}_{n-1}) - \mathcal{A}_h(\vec{x}^i_n)
\]
of the weak residual $r_n$ after the $j$th Newton iteration reaches $\rho > \rho_{\max}$. We use $\rho_{\max}=0.1$ or $\rho_{\max}=0.04$ in our experiments and stop when the residual satisfies $\vec{r}_n^{(j)}/\vec{r}_n^{(0)} < \num{1e-8}$ which is almost always reached in less than $10$ iterations.
The linear systems in each Newton step are solved with the GMRES method~\cite{kelleyIterativeMethodsLinear1995b} preconditioned with a V-cycle geometric multigrid~\cite{beckerMultigridTechniquesFinite2000}. Internally, a shared-memory parallel Vanka smoother is used, see~\cite{failerNewtonMultigridFramework2020} for more details on the efficient implementation.

%%%%%%%%%%%%%%%%%%%%%%%%%%%%%%%%%%%%%%%%%%%%%%%%%%%%%%%%%%%%%%%%%
\section{Deep neural network multigrid method}\label{sec:dnn-mg}
Given a numerical solver, e.g. as described above, the DNN-MG method uses a hierarchy of meshes $\Omega_h^0, \dots, \Omega_h^{L+J}$ that are obtained by successive refinement of the base mesh $\Omega_h^0$. 
We call $L\in \N_0$ the coarse level, $J \in \N$ the jump level and $L+J$ the fine level. On each level $l=0,\dots,L+J$ we define finite element spaces $\vec{V}_h^l$ and $\vec{L}_h^l$ as introduced in \cref{sec:discretization}. We can then use the numerical solver to cheaply compute a coarse solution on level $L$. The deep neural network is then used to refine it to a fine level $L+J$, which then feeds back to the coarse scales through the right hand side of \eqref{eq:nonlinearshort} during the time-stepping.

\begin{algorithm}[t]
    %\LinesNumbered%
    \SetFuncSty{textbf} \SetCommentSty{textsf} 
    \SetKwInOut{Input}{input} %
    \SetKwInOut{Output}{output} %
    \SetKwProg{Function}{function}{}{end} %
    \SetKwFunction{fnSolve}{Solve}%
    \SetKwFunction{rhs}{Rhs}%
    
    \Input{initial value $\vec{x}_0^{L+J}$, forcings $\vec{f}_0^{L+J}, \dots, \vec{f}_N^{L+J}$}
    \Output{solution $\nnsol{\vec{x}}_1^{L+J}, \dots, \nnsol{\vec{x}}_N^{L+J}$}
    \BlankLine
    $\nnsol{\vec{x}}_0^{L+J} = \vec{x}_0^{L+J}$ \tcp*{Initialize solution}%
    \label{alg:init}
    \For{$n=1, \dots, N$}{%
      $\vec{b}_{n}^{L+J} = \mathcal{F}_h^{L+J}(\nnsol{\vec{x}}_{n-1}^{L+J},\vec{f}^{L+J}_{n-1},\vec{f}^{L+J}_{n})$\tcp*{Evaluate the right hand side on level $L+J$}
      \label{alg:rhs}%
      $\vec{b}_{n}^{L} = \mathcal{R}(\vec{b}_{n}^{L+J})$\tcp*{Restrict rhs to level $L$}
      \label{alg:restrictrhs}%
      $\vec{x}_n^L$ = \fnSolve{$\vec{x}_{n-1}^L, \vec{b}_{n}^L$} \tcp*{Solve on coarse level}%
      \label{alg:coarsesolve}
      $\intsol{\vec{x}}_n^{L+J} = \mathcal{P}(\vec x_n^{L}) $\tcp*{Prolongate solution to level $L+J$}%
      \label{alg:prolong}
      $\vec{r}_n^{L+J}=\mathcal{F}_h^{L+J}(\intsol{\vec{x}}_n^{L+J}, \vec{f}_{n-1}^{L+J}, \vec{f}_{n}^{L+J})-{\mathcal A}_{h}^{L+J}(\intsol{\vec{x}}_n^{L+J})$ \tcp*{Compute residual on level $L+J$}
      \label{alg:residual}%
      $\nnsol{\vec{d}}_n^{L+J} =
        \mathcal{L}\circ\mathcal{N}(\mathcal{Q}(\tilde{\vec x}_n^{L+J}),\,\mathcal{Q}
        (\vec r_n^{L+J}),\,\Omega_{h}^{L+J})$\tcp*{Predict defect}
      \label{alg:predict}%
      $\nnsol{\vec{x}}_n^{L+J} = \intsol{\vec x}_n^{L+J} + \vec{\nnsol{d}}_n^{L+J}$\tcp*{Correct solution}
      \label{alg:correction}%
     % $\vec{b}_{n+1}^{L+J} = \mathcal{F}_h^{L+J}(\nnsol{\vec{x}}_n^{L+1},\vec{f}^{L+J}_n,\vec{f}^{L+J}_{n+1})$\tcp*{Right hand side for next time-step}
    %  \label{alg:rhs}%
     % $\vec{b}_{n+1}^{L} = \mathcal{R}(\vec{b}_{n+1}^{L+J})$\tcp*{Restriction of rhs
     %   to level $L$}
     % \label{alg:restrictrhs}%
    }%
  \caption{Time-stepping with DNN-MG.}
  \label{alg:dnnmg}
\end{algorithm}
The coupling of the numerical solver and neural network are described in \cref{alg:dnnmg}. When time-stepping with DNN-MG, we maintain a solution $\nnsol{\vec{x}}_n^{L+J}$ on the fine level $L+J$ (Line~\ref{alg:init}). As this is the ground truth that should steer the simulation, we evaluate the right hand side $\vec{b}_n^{L+J}$ on the fine level (Line~\ref{alg:rhs}). Afterward, the $L_2$-projection $\mathcal{R}: L_h^{L+J} \times V_h^{L+J} \to L_h^{L} \times V_h^{L}$ is used to restrict $\vec{b}_n^{L+J}$ to the coarse level (Line~\ref{alg:restrictrhs}), where we solve numerically (Line~\ref{alg:coarsesolve}). Then, the prolongation $\mathcal{P}:L_h^{L} \times V_h^{L} \to \mathcal{R}: L_h^{L+J} \times V_h^{L+J} $ transfers the solution $\vec{x}_n^L$ back to the fine level through interpolation (Line~\ref{alg:prolong}). The residual of this provisional solution $\intsol{\vec{x}}_n^{L+J}$ (Line~\ref{alg:residual}) serves as one of the inputs to the neural network that predicts a defect $\nnsol{\vec{d}}_n^{L+J}$ (Line~\ref{alg:predict}). The defect is finally added to the interpolated result to correct the solution for the current step (Line~\ref{alg:correction}).

\begin{figure}
    \centering
    \ifbuild
\else
	\includegraphics{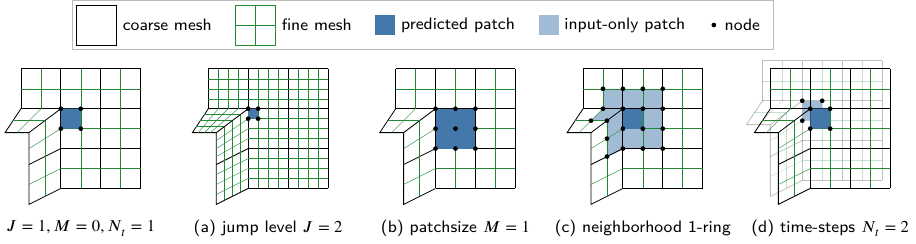}
	\fi
    \caption{Design choices regarding patches, which are the basic units for the prediction. For increased clarity, the indicated degrees of freedom (nodes) are for first order elements, although we always use quadratic elements  in our numerical experiments. The jump level (a) determines the number of mesh subdivisions done starting from the coarse grid for the numerical solution to reach the final degrees of freedom. A patch then corresponds to a cell of the fine mesh. As shown in panel (b), the patch size can be increased by combining cells, reversing the subdivision, while maintaining the degrees of freedom from the fine mesh. Another way to enlarge the receptive field of the neural network is to include multiple input patches from a local neighborhood, as shown in panel (c). Unlike larger patches, this extension does not need to conform to the hierarchical mesh structure, but it requires a suitable neural network architecture to deal with the variable number of inputs. The addition of patches from previous time-steps (d), on the other hand, is always possible.}
    \label{fig:patches}
\end{figure}
The neural network operates on local patches $P^M \subset \Omega^{L+J}_h$, a collection of cells $K \in \Omega^{L+J}_h$. The patch size $0 \leq M < L+J$ determines how many levels of refined cells are combined while keeping the degrees of freedom from the fine level. This is illustrated in \cref{fig:patches}, where $P^0$ corresponds to a cell of the fine mesh, while $P^J$ corresponds to a cell of the coarse mesh. Each patch consisting of quadratic quadrilaterals thus contains 
\begin{equation}\label{eq:nodespatch}
N_P^M = (2^{M+1} + 1)^d
\end{equation}
nodes, where  $d\in\{2,3\}$ is the dimension of the problem. To obtain the input for the neural network, 
\[
\mathcal{Q}: N^{L+J}_{\text{node}} \times N_\text{comp} \to N^{L+J-M}_{\Omega} \times N_P^M N_\text{comp}
\]
gathers the associated $N_{\text{comp}}$ degrees of freedom at each of the $N^{L+J}_{\text{node}}$ nodes from a global vector and produces one local vector for each of the $N_{\Omega_h}^{L+J-M}$ patches corresponding to cells on the mesh $\Omega_h^{L+J-M}$. The local vectors are the inputs for the neural network 
\[
\mathcal{N}: (\dots) \to N^{L+J-M}_{\text{node}} \times N_\text{comp}
\] 
which receives the interpolated solution $\mathcal{Q}(\tilde{\vec x}_n^{L+J})$, the residual $\mathcal{Q}(\vec r_n^{L+J})$ and the cell geometry $\Omega_{h}^{L+J}$ as a batch and predicts the velocity defect. The concrete shape of the inputs depends on a number of parameters and further processing is applied, the details of which are discussed in \cref{sec:io}. Afterward, the neural network is evaluated on each of the $N^{L+J-M}_{\text{node}}$ patches in parallel. The resulting local corrections have to be scattered again to a global vector by
\[
\mathcal{L}: N^{L+J-M}_{\Omega} \times  N_P^M N_{\text{comp}} \to N^{L+J}_{\text{node}} \times N_\text{comp}.
\] 
Care is needed during the scattering operation because the degrees of freedom at the edges of patches appear multiple times in different local vectors. 
This is resolved by summing the contributions from multiple vectors and rescaling by the reciprocal of their number of appearances. For a more rigorous description of the handling of the patches see~\cite{Margenberg2024}.

%%%%%%%%%%%%%%%%%%%%%%%%%%%%%%%%%%%%%%%%%%%%%%%%%%%%%%%%%%%%%%%%%
\subsection{Inputs and outputs}\label{sec:io}
The patch structure of the neural network input in DNN-MG gives us two important parameters, shown in \cref{fig:patches}. 
These are crucial for the supervised learning task of the neural network.  

The jump level $J$ must be chosen based on the desired accuracy of the solution. If the predicted resolution is too low to properly resolve important features, the correction will not be effective. 
%For example, $J=2$ levels were needed for Navier-Stokes in 3D in \cite{Margenberg2024}. 
On the other hand, a higher jump level can also degrade the prediction accuracy because it shrinks the area covered by a patch. The neural network thus has comparably less information available from the coarse solution on each patch. In previous work, $J=2$ proved effective ~\cite{Margenberg2024, Margenberg2025}, but given the added cost of data generation for each level, we use $J=1$ for most of our experiments. However, a discussion regarding the effect of a reduced coarse level is provided in \cref{sec:generalization}.

A remedy for the lack of spatial context when solving a non-local problem can be to increase the patch size $M$. By setting $M=J$, the neural network will always see the complete coarse cell and if the mesh has more than two levels, a patch size $M>J$ is also possible. 
%In fact, a-priori error estimates in \cite{Kapustsin2023} indicate that  
Our mathematical analysis for simplified problems~\cite{Kapustsin2023} even suggests that the ratio $M/J$ must grow for decreasing mesh sizes in order to be able to guarantee robustness. It is, however, currently unclear if this is a weakness of the analysis as our numerical experiments show robustness even when $M \leq J$.
A larger patch size makes the task more difficult from a deep learning perspective. Since the number of nodes included in a patch increases according to \eqref{eq:nodespatch}, the input and output spaces of the neural network grow by the same factor. Thus, one expects that more training samples are needed to achieve a similar generalization error~\cite{golestaneh2025how}, while at the same time, the number of samples derived from a simulation is reduced by a factor of $4$ for each additional patch level. We compare different patch sizes in practice in~\cref{sec:generalization}. Unless explicitly stated, we use $M=0$ in the experiments.

Another way to increase the receptive field of the neural network is to use multiple patches as inputs while still predicting a single patch. Both spatial and temporal extensions as in \cref{fig:patches} (c) and (d) are possible. In contrast to an increased patch size, this approach does not reduce the number of training samples and in the case of (c), we also ensure that even the degrees of freedoms at the edges of the predicted patch have context available in all directions. A fixed number of input patches can be used with any neural network while a varying number requires specialized architectures such as graph neural networks or transformers, cf. \cref{sec:nndesign}. 

\subsection{Patch inputs}
We now turn to the concrete data prepared for each input patch $P^M$. In line with our previous work~\cite{margenbergDeepNeuralNetworks2021,Margenberg2024}, we provide three types of inputs:
\begin{enumerate}
    \item the nonlinear residuals $\vec{r}^{L+J}_{n,o} \in \R^{N_{\text{comp}}}$ of the prolongated coarse solution~(\cref{alg:dnnmg} Line~\ref{alg:residual}) for each node $o\in P^M$,
    \item the prolongated coarse solution $\vec{x}^{L+J}_{n,o} \in \R^{N_{\text{comp}}}$ (\cref{alg:dnnmg} Line~\ref{alg:prolong}) for each $o \in P^M$,
    \item the geometry of the patch $P^M$ in the form of the $2^d$ relative corner positions $\vec{g}_c^* \in \R^d, c \in K$ with
    \[
    \vec{g}_c^* = \vec{g}_c - \frac{1}{2^d}\sum_{u\in K} \vec{g}_u,
    \]
    where $\vec{g} \in \R^{N^{L+J-M}_\Omega \times d}$ are the vertices of the mesh $\Omega_h^{L+J-M}$ and $K \in \Omega_h^{L+J-M}$ is the cell corresponding to the whole of $P^M$.
\end{enumerate}
A notable change compared to our previous work is that we normalize all inputs and outputs. Scale and shift values are estimated from the full training data to get each component close to zero mean and unit variance. This is standard practice in deep learning, since it ensures that the neural network is sensitive to all components of the input and since it stabilizes training. In particular, the solution is naturally of order $1$ but the residual is more than $3$ orders of magnitude smaller, so that previously it had potentially only a limited effect. The new normalization also influences the representation of the geometry. With it, the relative vertex positions $\vec{g}^*_c$ perform better than the combination of edge lengths, aspect ratios and angles used in \cite{MargenbergRichter2020}.

We investigated one more modification to the inputs. When at least one previous time-step is provided, i.e. $N_t>1$, it is possible to add the neural network corrections back to the inputs for successive steps. Instead of providing the prolongated coarse solution $\vec{\intsol{x}}^{L+J}_{n-1}$, we can insert $\vec{\nnsol{x}}^{L+J}_{n-1}$. In essence, this turns the task of the neural network from downscaling into time-stepping, since at least an approximation of the full fine state on $\Omega_h^{L+J}$ is available. When the neural network also receives patches in the 1-ring neighborhood, it can effectively solve the transport problem. Indeed, when we combine temporal and spatial multi-patch inputs with fine velocity inputs for past steps, we see improvements in the final validation loss, in some cases by more than a factor of $2$. However, this improvement does not translate into increased accuracy in the hybrid simulation with the finite element solver. The direct feedback from the neural network's own corrections quickly causes errors to accumulate, which destabilizes the simulation. In our experiments, the simulation quickly developed nonphysical features and usually exploded after a few hundred steps. The implementation of fine solution inputs and our limited success to stabilize this setup are further discussed in \cref{appendix:fine-vel-inputs}.

\subsection{Training}\label{sec:training}
The neural network is trained on data from reference simulations. The samples are generated by running the solver in much the same way as the hybrid simulation described in \cref{alg:dnnmg}. The solution $\intsol{\vec{x}}^L_n$ for the current time-step is first computed on the coarse level, prolongated to the fine level as $\intsol{x}^{L+J}_n$ and the inputs for the neural network are prepared. Then, instead of evaluating a neural network, the numerical solver is run on the fine level to compute the reference solution $\refsol{\vec{x}}^{L+J}_n$. From this solution, the defect
\[
\refsol{\vec{d}}_n^{L+J} = \refsol{\vec{x}}^{L+J}_n - \intsol{\vec{x}}^L_n
\]
is computed, which is the target in our supervised training. Afterward the input-target pairs for the whole domain are written out as training data.

During the training, we built mini-batches from random samples of multiple simulations. One sample corresponds to the state of a single patch in a single time-step so that mini-batches combine samples from different times and even different cases. As one epoch we consider training on all samples of the training set once. As validation set, we use simulation results from cases with unique geometry that are not part of the training set. After each epoch, the mean square error is evaluated on the validation data to monitor the changes and to select the best neural network in the end. The specific cases used and the resulting dataset sizes are reported as part of the experiment descriptions \cref{sec:noise-aug} and \cref{appendix:experiment-details}. Further details on the loss function and the optimizer are dependent on the neural network architecture and will be described in \cref{sec:nndesign}.

%%%%%%%%%%%%%%%%%%%%%%%%%%%%%%%%%%%%%%%%%%%%%%%%%%%%%%%%%%%%%%%%%
\subsection{Neural network design}\label{sec:nndesign}
The DNN-MG method can be realized with a wide range of neural network architectures. Since the same neural network is used on every patch, approximate translational invariance, similar to that of convolutional neural networks, is already built into the method. 
\begin{table}
\begin{tabular}{lcccc|cccc}
\toprule
               & \multicolumn{4}{c}{$d=2$}     & \multicolumn{4}{c}{$d=3$}     \\
               & $M=0$ & $M=1$ & $M=2$ & $M=3$ & $M=0$ & $M=1$ & $M=2$ & $M=3$ \\
\midrule
$N_\text{in}$  & 62    & 158   & 494   & 1742  & 240   & 1024  & 5856  & 39328 \\
$N_\text{out}$ & 18    & 50    & 162   & 578   & 81    & 375   & 2187  & 14739 \\
\bottomrule
\end{tabular}
    \caption{Input and output sizes of the neural network for different patch sizes $M$ and dimensions $d$. Each vector has $N_\text{comp}=d+1$ (velocity+pressure) components per node and $N^M_P$ nodes per patch (\cref{eq:nodespatch}), which taken together with the geometry results in $N_\text{in} = 2 N^M_P N_\text{comp} + 2^dd$ inputs. The outputs are just the defect of the velocity components, totaling $N_\text{out} = N^M_P d$. The jump level $J$ does not impact these sizes directly but it may influence the choice of $M$.}
    \label{tab:in-out-sizes}
\end{table}
This also leads to small input and output dimensions, see~\cref{tab:in-out-sizes}. In this work, we mostly use $J=1, M=0$ in 2D, giving $N_\text{in}=62$ and $N_\text{out}=18$. The patch wise approach thus enables a compact neural network design with a small number of trainable parameters and comparably efficient training. In previous work, we successfully applied a simple multilayer perceptron (MLP) to the 3D Navier-Stokes equations \cite{Margenberg2024} with $J=2, M=1$. Exploring the use of larger $M$, where dense layers become impractical, we leave for future work.

While simply scaling up the size of a feedforward neural network can improve the accuracy of DNN-MG, other architectures with stronger inductive biases could further improve the results by efficiently utilizing multiple input patches. We therefore explored three different neural network designs with different capabilities in our experiments. 
Their design and different properties are summarized below and described in more detail in \cref{appendix:nn-details}.

\subsubsection{Multilayer perceptron}
As baseline we use a simple multilayer perceptron (MLP). An $\text{MLP}: \R^{N_\text{in}} \to \R^{N_\text{out}}$ is a sequence of linear maps with learnable weights, each followed by a fixed nonlinear activation function that is applied elementwise. After each linear map we also apply a layer norm~\cite{Ba2016}. In the hidden layers, where the input and output dimensions match, we insert skip connections that add the input directly to the output and which stabilizes training.
 
We use a layer norm (instead of a batch norm as in~\cite{Margenberg2024}) for its simplicity, since training and inference are inherently aligned, and to be consistent with the other architectures where batch normalization is not always applicable. In general, normalization accelerates training and improves the generalizability of a neural network; see~\cite{Huang2023} for an extensive survey. In DNN-MG, normalization is important to make the hybrid simulation with the MLP stable. As activation function we use tanh, as it performed better in our experiments than commonly used rectified activations such as swish~\cite{Ramachandran2017}.

The loss function used for training is the standard mean square error (MSE). A constant learning rate of \num{0.0005}, a batch size of 512, and the AdamW optimizer~\cite{loshchilovDecoupledWeightDecay2018} with a weight decay of $\num{5e-05}$ are used. AdamW turned out to be essential for stability during inference. In contrast to regular Adam~\cite{kingmaAdamMethodStochastic2015}, weight decay (or Tikhonov regularization) is decoupled in AdamW.

%\blue{The hyperparameters where tuned over multiple small scale grid searches, varying just one or two parameters at a time until good convergence was observed. We always did multiple runs with different random seeds and found these parameters to be robust across neural network sizes, which ranged from \num{1e+5} to \num{2e+6} weights in our experiments.} 

% batchsize 256 1024 also fine
% wd 1e-4 fine, 4e-4 too much
% smaller lr (0.0001) is bad

\subsubsection{Recurrent neural network}
Recurrent neural networks (RNN) maintain a hidden state that is based on previous time-steps. It serves as memory that is used together with the current step's input to make a prediction. In theory, the memory allows the neural network to retain information on the predicted fine solution, ensuring temporal consistency of the simulation. 
%Otherwise, the high-resolution solution features have to be reconstructed from the coarse solution in every step. 
This was the rationale for using a recurrent neural network in the original DNN-MG in~\cite{margenbergNeuralNetworkMultigrid2022}. The current design is an evolution of this and also based on Gated Recurrent Units (GRUs)~\cite{Cho2014}. Our $\text{RNN}: \R^{N_\text{in}} \times \R^{(L-1) \times N_\text{hid}} \to \R^{N_\text{out}} \times \R^{(L-1) \times N_\text{hid}}$ follows the same blueprint as the MLP. The linear maps are replaced by standard GRUs without bias. Layer normalization is applied between GRUs and skip connections allow the GRUs to learn an update instead of the full new state. The final layer is just a linear projection that maps the latent vector to the output space. 

For training, the same setup is used as for the MLPs, with the notable difference that BPTT is needed to learn RNN-type architectures with memory. Mini-batches are constructed from random cells with a fixed sequence length $w\in \N$ and the forward passes are rolled out for the entire sequence before the gradient update. Since inputs are still taken from the training data, the gradients can only flow through the hidden state updates. Even in this restrictive form, we note that BPTT can aid with simulation stability, as the RNNs can be trained effectively with simple Adam, in contrast to the MLPs where AdamW is required. Furthermore, the training benefits from long sequences with visible improvements in the validation loss until around $w=128$, the setting we use for subsequent experiments.

\subsubsection{Transformer}\label{sec:transformers}
Transformers~\cite{NIPS2017_3f5ee243} were developed mainly for faster sequence processing via increased parallelism compared to recurrent neural networks. Their weaker inductive bias makes them well-suited also for spatial data~\cite{Palanisamy2025}. Transformers operate on a sequence of vectors called tokens which, in our case, correspond to patches. In addition to previous time-steps of the current patch, we collect all current time inputs for patches in a 1-ring neighborhood, i.e., every patch that shares at least one vertex on the coarse mesh with the current one. We thus retain the benefits of the local approach, but also enlarge the receptive field of the neural network such that it has information from a neighborhood around every node for which it makes a correction. 

Our transformer $\text{Transf}: \R^{N_{\text{seq}} \times N_\text{in}} \to \R^{N_\text{out}}$ is based on a standard BERT style transformer~\cite{devlin-etal-2019-bert}. The BERT layers consist of alternating MLPs, applied to each token independently and self-attention, where the key operation are scaled dot-products between pairs of tokens. Our embedding that turns input patches into tokens is a learned linear transformation. A simple learned positional encoding is added to each token to differentiate between them. We also experimented with 2D positional encodings based on the actual positions relative to the target patch but found that they performed worse. The most important feature of the encoding seems to be that it identifies the token corresponding to the target patch. In order to force the neural network to collect relevant information from every patch, another learned token, commonly called the \texttt{[CLS]} token, is appended to the sequence and this token is later used to make the prediction.

Similar to \cite{Lessig2023}, the actual prediction is made by an an ensemble of $N_\text{ens} \in \N $ shallow MLPs and the final result is taken to be their mean. This $\text{Tail}: \R^{N_\text{seq} + 1} \times \R^{N_\text{tok}} \to \R^{N_\text{out}}$  selects the \texttt{[CLS]} token $\vec{k}_\text{cls} \in \R^{N_\text{tok}}$ and computes
\[
\text{Tail}(\vec{x}) = \frac{1}{N_\text{ens}} \sum_{j=1}^{N_\text{ens}} \text{MLP}_j\left(\vec{k}_\text{cls}\right).
\]
We use the standard GELU activations in both the transformer layers and the tail MLPs. In all the presented results we use an ensemble of size $N_\text{ens} = 4$.

\paragraph{Probabilistic training}
The tail ensemble makes the neural network more robust by allowing us to recast the training as probabilistic in order to account for the inherent uncertainty in the prediction. While the ground truth $\refsol{d} \in \R^{N_\text{out}}$ is just a single, deterministic trajectory, the neural network only receives local, lower resolution information to solve the non-local Navier-Stokes equation. Each prediction $\nnsol{d}_j \in \R^{N_\text{out}}$ of an ensemble member is hence considered as a sample drawn from the distribution of possible solutions and $\text{Transf}$ models the whole distribution. Treating the components $\nnsol{d}_j^i, i=1,\dots N_\text{out}$ as independent and normal distributed according to $\nnsol{d}_j^i \sim N(\nnsol{\mu}_i,\nnsol{\sigma}_i^2)$, we can compute the parameters
\[
\nnsol{\mu}_i = \frac{1}{N_\text{ens}} \sum_{j=1}^{N_\text{ens}} \nnsol{d}_j^i, \qquad 
\nnsol{\sigma}_i^2 = \frac{1}{N_\text{ens}} \sum_{j=1}^{N_\text{ens}} \left(\nnsol{\mu}_i - \nnsol{d}_j^i\right)^2.
\]
To train this distribution, we use the loss function $\mathcal{L}: \R^{N_\text{out}} \times \R^{N_\text{out}} \times \R^{N_\text{out}} \to \R$ with
\begin{equation}\label{eq:proploss}
\mathcal{L}(\vec{\refsol{d}}, \vec{\nnsol{\mu}}, \vec{\nnsol{\sigma}}) = \frac{1}{N_\text{out}} \sum_{i=1}^{N_\text{out}} \left( \frac{1}{2} (\refsol{d}_i - \nnsol{\mu}_i)^2 + \frac{1}{2}\left(\left( 1 - G_{\nnsol{\mu}_i,\nnsol{\sigma}_i}(\refsol{d})\right)^2 + \nnsol{\sigma}_i\right) \right),
\end{equation}
where $G_{\mu_i,\sigma_i}$ is unnormalized Gaussian probability density function
\[
G_{\mu_i,\sigma_i}(\vec{x}) = e^{-\frac{(\vec{x} - \mu_i)^2}{2 \sigma_i^2}}.
\] 
The first term is just the MSE of the whole ensemble. The second term, which estimates the distance between the predicted distribution and the unknown reference distribution from individual training samples $\bar{d}$, is taken from \cite{Lessig2023}. However, we require a stronger regularization term $\nnsol{\sigma}_i$ for stable training. We also tried the more common Gaussian negative log likelihood loss~\cite{Nix1994} but \cref{eq:proploss} gave better results.

Further parameters are important for robust training. Simply using the same hyperparameters as for the MLPs, we found that the transformers tend to overfit quickly, a problem that we did not encounter with the other architectures. After just a few epochs there would be a sharp increase in both loss and the standard deviation of the ensemble members on the validation sets. However, with some adjustments it is possible to achieve good convergence on both training and validation sets consistently.
Following the standard for transformers we do a warmup, starting with a low learning rate of $r_\text{min} = \num{2e-05}$ and apply the linear learning rate schedule 
\[
r(c) = \begin{cases} 
	r_\text{min} + \frac{c}{c_\text{up}} (r_\text{max} - r_\text{min}) & c \leq c_\text{up} \\
	r_\text{min} + (1 - \frac{c-c_\text{up}}{c_\text{down}}) (r_\text{max} - r_\text{min}) & c > c_\text{up}
 \end{cases}
\]
for epochs $c = 0, \dots, 63$. The learning rate increases to $r_\text{max} = \num{1.8e-4}$ at $c_\text{up}=14$ and slowly decays again until $c_\text{down} = 120$. Furthermore, compared to MLP training, we use a larger batch size of $2048$ and a stronger weight decay of \num{0.0001}.

%%%%%%%%%%%%%%%%%%%%%%%%%%%%%%%%%%%%%%%%%%%%%%%%%%%%%%%%%%%%%%%%%%%%%%%%%%%%
\section{Test cases and metrics}\label{sec:cases}
\begin{figure}
    \centering
\ifbuild
\else
\includegraphics{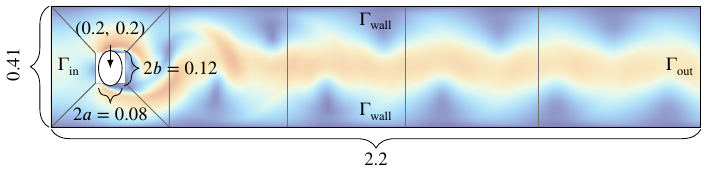}
\fi
\ifbuild
\else
\includegraphics{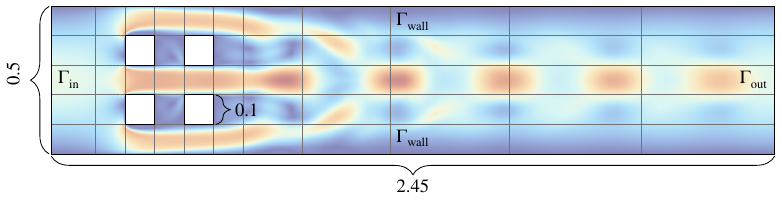}
\fi
    \caption{Geometry of two example cases (ro1, sq4). Each channel has a parabolic inflow profile $\Gamma_{\textrm{in}}$,
  do-nothing boundary conditions at the outflow boundary $\Gamma_{\textrm{out}}$ and no-slip conditions on the walls
  $\Gamma_{\textrm{wall}}$ and the obstacles. The thin lines indicate the base mesh which is refined uniformly 5 times for round obstacles, conforming to the elliptical shape, and 4 times for square obstacles to create the coarse grid for the numerical simulations. In the background, a velocity magnitude snapshot of a fully developed fine simulation is shown.}
  \label{fig:meshes}
\end{figure}
In our numerical experiments, we work with channel flows of different geometries. In \cref{fig:meshes} we show the two main types of meshes from which the cases are built. Elliptical obstacles are created through conforming refinement~\cite{Gascoigne3D}, square obstacles are simply holes in the otherwise rectangular mesh. Variations are created by arranging multiple obstacles with the same local mesh in a grid and extending the channel as needed in width and length. For round obstacles, we furthermore vary the ellipse radii $a,b \in [0.35, 0.65]$ and the exact y-position within the square by $\{-0.02, 0.0, 0.02\}$. The resulting cases will be abbreviated as ``ro$N$'' for round obstacles and ``sq$N$'' for square obstacles, where $N$ is the number of obstacles. Detailed descriptions for each mesh are given in \cref{tab:mesh-details}.

The flow is always driven by a Dirichlet profile $v=v^D$ at the left boundary $\Gamma_{\mathrm{in}}$, given by
\[
 \begin{aligned}
 % v_\text{avg}\frac{6y(H-y)}{H^2}\omega(t)
v^D(x,y,t) &= \frac{y(H-y)}{H^2}\omega(t) \ \text{  on  } \ \Gamma_{in} :=
0\times [0,\,H],
\\[4pt]
\; \omega(t) &= \begin{cases} \frac{1}{2}- \frac{1}{2}\cos( 2\pi t), & t\le \frac{1}{2} \\ 1 & t>\frac{1}{2},
\end{cases}
\end{aligned}
\]
where $H$ is the height of the channel and the function $\omega(t)$ regularizes the start-up phase of the flow. On the wall boundary $\Gamma_{\text{wall}}$ and on the obstacle we prescribe no-slip boundary conditions $\vec{v}=0$. On the outflow
boundary $\Gamma_{\mathrm{out}}$ we use a do-nothing outflow condition~\cite{HeywoodRannacherTurek1992}.

The flow that results with the above setup has a Reynolds number of around $\Rey\approx100$, leading to vortices being created by the round obstacles. An arrangement of square obstacles, on the other hand, develops a stationary solution in the front part of the channel and vortices only occur Afterward in the open channel. The re-entrant corners of the square obstacles limit the regularity of the solution and hence we expect reduced convergence rates on finer meshes. 

To properly evaluate DNN-MG and its ability to generalize to flow conditions not seen during training, it is necessary to look at different cases and to use a variety of metrics. In practice, the accuracy of the neural network prediction can vary wildly between seemingly similar scenarios. Since the neural network introduces an uncontrolled error, conventional theory of convergence based on the grid spacing is not applicable.  

A particular useful case is the well-established ``laminar ﬂow around a cylinder''~\cite{SchaeferTurek1996}. The unsteady 2D-2 case (ro1 in our naming scheme) settles into a stable time-periodic flow. The solution can thus be well characterized by the drag and lift functionals
\begin{align}
J_d(\vec{v}_n,p_n) = \int_\Gamma \Big(\frac{1}{\Rey} \nabla \vec{v}_n - p_nI\Big) \vec{n} \cdot \vec{e}_1 \,\text{d}s,\label{eq:drag}\\
J_l(\vec{v}_n,p_n) = \int_\Gamma \Big(\frac{1}{\Rey} \nabla \vec{v}_n - p_nI\Big) \vec{n} \cdot \vec{e}_2\,\text{d}s,\label{eq:lift}
\end{align}
where $\vec{e}_1=(1,0)^T, \vec{e}_2=(0,1)^T$ and $\vec{n} \in \R^2$ is the surface normal. A direct comparison of individual time-steps from different simulations, e.g. via the $L_1$ norm
\begin{align}
e_{J_d}(\refsol{\vec{v}}_n, \refsol{p}_n,\nnsol{\vec{v}}_n, \nnsol{\vec{v}}_n) &= |J_d(\refsol{\vec{v}}_n, \refsol{p}_n) - J_d(\nnsol{\vec{v}}_n, \nnsol{\vec{v}}_n) |,\label{eq:drag-error}\\
e_{J_l}(\refsol{\vec{v}}_n, \refsol{p}_n,\nnsol{\vec{v}}_n, \nnsol{\vec{v}}_n) &= |J_l(\refsol{\vec{v}}_n, \refsol{p}_n) - J_l(\nnsol{\vec{v}}_n, \nnsol{\vec{v}}_n) |, \label{eq:lift-error}
\end{align}
only makes sense at the beginning of the simulation because the frequency of the oscillation increases with the mesh resolution~\cite{MargenbergRichter2020} and even a small phase change causes large differences in the $L_1$ norm. 
To compare the fully developed flow, we therefore use aggregate values over a time interval $[t_{n_0}, t_{n_0+k}]$ and additionally estimate the frequency of the solution from the lift functional. A good approximation can be computed with only a few periods by first finding all steps $n_i, i=1,\dots,m$ within the interval where $J_l(v_{n_i},p_{n_i}) \leq 0$ and $J_l(v_{n_i+1},p_{n_i+1}) > 0$. Then, we use a linear approximation through $n_i$ and $n_{i+1}$ to estimate the roots $\hat{t}_{n_i}$ of $J_l$, from which we compute the frequency
\begin{equation}\label{eq:freq}
    f(v, n_0, k) = \frac{m-1}{\hat{t}_{n_m} - \hat{t}_{n_1}}.
\end{equation}

The functionals allow us to quantify the accuracy of the complicated flow close to the obstacle. However, in order to evaluate the whole velocity field, especially for cases with more chaotic solutions, other metrics are needed. Another value that can be computed efficiently during the simulation is the divergence of the velocity vector field
\begin{equation}\label{eq:div}
    J_\text{div}(\vec{v}_n) = \int_\Omega (\nabla \cdot \vec{v}_n)^2 \,\text{d}\vec{x}.
\end{equation}
Since we are solving the incompressible Navier-Stokes equation, the divergence should approximately vanish. The stabilization term in \cref{eq:4cranknicholson3} introduces some divergence but the violation of divergence-freedom should reduce by a factor of 4 for each additional mesh level. Since there only small fluctuations of $J_\text{div}$ within one simulation, the divergence can be easily compared between reference and DNN-MG simulations. 

A more direct comparison of the velocity fields in the presence of a phase shift is possible by computing the mean velocity over an interval $[t_{n_0}, t_{n_0+k}]$ and taking the $L_2$ norm of the difference
\begin{equation}\label{eq:mean-vel-error}
    e_{\bar{v}}(\refsol{\vec{v}}, \nnsol{\vec{v}}, n_0, k) = \frac{\left\lVert \frac{1}{k}\sum_{n = n_0}^{n_0 + k} \refsol{\vec{v}}_n - \frac{1}{k}\sum_{n = n_0}^{n_0 + k} \nnsol{\vec{v}}_n \right\rVert_2}{ \left\lVert \frac{1}{k} \sum_{n = n_0}^{n_0 + k} \refsol{\vec{v}}_n \right\rVert_2}.
\end{equation}
The values $n_0$ and $k$ have to be carefully chosen for $e_{\bar{v}}$ to be meaningful. To ensure that the variation in time does not play a role, we compute the time-averaged velocity $e_{\bar{v}}$ for increasing $k$ until we observe convergence.

A useful error metric that completely sidesteps phase shift issues is the local error. It can be computed for a single simulation as 
\begin{equation}\label{eq:trunc-error}
    \tau(\nnsol{\vec{v}}_{n-1}, \nnsol{\vec{v}}_n, \vec{f}_{n-1}, \vec{f}_n) = \frac{\lVert \mathcal{A}_h^{-1}(\mathcal{F}_h(\nnsol{\vec{v}}_{n-1}, \vec{f}_{n-1}, \vec{f}_n))_\vec{v} - \nnsol{\vec{v}}_n \rVert_2}{\lVert \mathcal{A}_h^{-1}(\mathcal{F}_h(\nnsol{\vec{v}}_{n-1}, \vec{f}_{n-1}, \vec{f}_n))_\vec{v} \rVert_2}.
\end{equation}
Effectively, we assume the previous time-step $\nnsol{\vec{v}}_{n-1}$ to be the correct state and compute a ground truth solution only for the current step. In practice, the computation of the reference is very similar to the training data generation procedure. First, we solve on the coarse grid as input to the neural network and compute the corrected solution $\nnsol{\vec{v}}_n$. Then we compute the reference solution
\[
\refsol{\vec{v}}_n = \mathcal{A}_h^{-1}(\mathcal{F}_h(\nnsol{\vec{v}}_{n-1}, \vec{f}_{n-1}, \vec{f}_n))_\vec{v}
\]
with the numerical solver on the fine grid using the same right hand side. In contrast to the data generation procedure, the fine reference solution $\refsol{\vec{v}}_n$ is discarded Afterward and neural network corrected solution $\nnsol{\vec{v}}_n$ is taken to compute the right hand side for the next step. There are two problems with \eqref{eq:trunc-error}, it is expensive to compute and it does not show systemic deviations over time that can make the flow nonphysical. 
%relative velocity error
%\begin{equation}
%    e_v(\refsol{v}_n, \nnsol{v}_n) = \frac{\lVert \refsol{v}_n - \nnsol{v}_n \rVert_2}{\lVert\refsol{v}_n \rVert_2}
%\end{equation}

%%%%%%%%%%%%%%%%%%%%%%%%%%%%%%%%%%%%%%%%%%%%%%%%%%%%%%%%%%%%%%%%%%%%%%%%%%%%
\section{Data augmentation and stability}\label{sec:stability}
With the recipe described above, DNN-MG is generally stable in the usual sense of a numerical method. The perturbations introduced by the neural network do not disrupt the convergence of the numerical solver and even over many time-steps there's no excessive error growth. However, instability can manifest in a different way in autoregressive time-stepping with a neural network. Although a sufficiently regularized neural network is unlikely to lead to an exponential feedback loop as conventional simulation do, the small errors that are introduced at each step can cause a shift over time. This is what we observe in \cref{fig:solution-nn-artifacts}, where nonphysical high frequency artifacts appear in the solution after a few steps. For the incompressible Navier-Stokes equations, the divergence of the solution is a good proxy for this error. Our finite element solutions are not divergence-free, but each mesh refinement should reduce the divergence by a factor of 4. The neural network learns this structure and initially reduces the divergence compared to the coarse solution. However, the divergence in \cref{fig:div-nn-artifacts} grows rapidly and eclipses that of the interpolated coarse solution after just a few steps. At this point, the simulation leaves the flow regime on which the neural network was trained and the network will produce a large error at each step. In the example, after the divergence peak around $t=\qty{1.2}{s}$, the simulation settles into a different stable equilibrium but one that has noticeable artifacts.

A DNN-MG solution with visible artifacts is not necessarily worse than the coarse solution by most metrics, as the effect on the macroscopic flow is usually limited; cf. \cref{fig:div-noisy-aug}. Nevertheless, this distributional shift constitutes a large source of errors which is not reduced by improving the single-step prediction that we train on. In fact, the problem can get worse with better neural networks as measured by the validation loss, since these are in general more sensitive to the inputs and their quality. In theory, it is possible to minimize the cumulative error directly with backpropagation through time (BPTT). For long roll-outs, the loss in later time-steps is increasingly dominated by the incorrect input from previous steps instead of the accuracy of the current prediction. The main problem in our case is that BPTT would require the computation of gradients throughout the whole numerical scheme. 
Making an existing code base like Gascoigne differentiable is difficult. One possible approach is to re-implement both the conventional solver and the neural network within an autograd framework such as jax~\cite{jax2018github} or PyTorch~\cite{Paszke_PyTorch_An_Imperative_2019}, as has been done, e.g., in \cite{Kochkov2024}. Another alternative is the use of a modified compiler like Enzyme~\cite{NEURIPS2020_9332c513} with the existing code. In addition to the considerable implementation effort, performance is of practical concern for both approaches. The gradient computation itself has a significant cost~\cite{Moses2022} and an implementation with an autograd framework can be much slower than a handwritten C++ implementation~\cite{Jendersie2025}.
Therefore, we look to alternative methods which do not require gradients in the solver, but instead work by augmenting the training data, to alleviate the issue of cumulative errors.
%This could be achieved by implementing both the conventional solver and the neural network within a framework such as jax~\cite{jax2018github} or PyTorch~\cite{Paszke_PyTorch_An_Imperative_2019}, as has been done, e.g., in \cite{Kochkov2024}. Alternatively, a modified compiler like Enzyme~\cite{NEURIPS2020_9332c513} could be used to make the existing code differentiable. However, in practical terms the performance would be a significant issue with either approach. Code written in such frameworks can be much slower than a handwritten C++ implementation~\cite{Jendersie2025} and the gradient computation itself has a considerable cost~\cite{Moses2022}. This is exaggerated by the complexity of the multigrid method and the fact that the implicit solver always requires the global solution. Therefore, we look to alternative methods which do not require gradients in the solver, but instead work by augmenting the training data, to alleviate the issue of cumulative errors. 

%%%%%%%%%%%%%%%%%%%%%%%%%%%%%%%%%%%%%%%%%%%%%%%%%%%%%%%%%%%%%%%%%
\subsection{Added noise}\label{sec:noise-aug}
A simple way to increase the variety of samples in the training set is to add noise to the simulation during the data generation. Perturbations to the fine grid velocity approximate the effect of the neural network in the later hybrid simulation. Thus, as long as the numerical solver still converges with the perturbations, it will steer the simulation back to a physical state. Noises samples can therefore teach the neural network to be robust to its own errors.

Our noise schedule adds random values distributed according to $N(0, \sigma^2)$ to each velocity component of the fine solution after the training data for the current step is written. Since it can take multiple steps for the simulation to reach the baseline divergence again after a perturbation, we only add noise with probability $0 \leq c^{-1} \leq 1$ per time-step. The resulting dataset thus still has samples very close to the unperturbed simulation and the random period ensures that no spurious patterns are introduced. This technique is easy to implement and cheap since the added cost during the data generation is negligible. In principle, one could combine unperturbed and noisy samples for a dataset, but since each requires a complete simulation, the compute time is likely better spend on more distinct cases. %we still get the same number of training samples and the added cost during the data generation is negligible.

\begin{figure}
    \centering
    \ifbuild
	\tikzsetnextfilename{noise_aug_ro1_div}
	\begin{tikzpicture}
		\begin{axis}[%
			legend pos=north east,
			legend columns = 2,
			xminorticks=true,
			xmin=-0.005,
			xmax=0.045,
			width=0.49\textwidth,
			xlabel={$\sigma$},
			ylabel={mean $J_\text{div}$},
			]
			\addplot+[mark=none, dashed] table[x=setups, y=mean_div_coarse]{functionals_gc1_noisy_data_reference.txt};
			\addplot+[mark=none, dashdotted] table[x=setups, y=mean_div_fine]{functionals_gc1_noisy_data_reference.txt};
			
			\addplot+[error bars/.cd, y dir=both, y explicit] table[x=noise_sigma,y=mean, y error=std] {functionals_gc1_noisy_data_pr1_mean_div_lines.txt};
			\addplot+[error bars/.cd, y dir=both, y explicit] table[x=noise_sigma,y=mean, y error=std] {functionals_gc1_noisy_data_pr8_mean_div_lines.txt};
			\addplot+[error bars/.cd, y dir=both, y explicit] table[x=noise_sigma,y=mean, y error=std] {functionals_gc1_noisy_data_pr16_mean_div_lines.txt};
			
			\legend{coarse, fine, $p=1$, $p=8$, $p=16$};
		\end{axis}
	\end{tikzpicture}
    \else
	\includegraphics{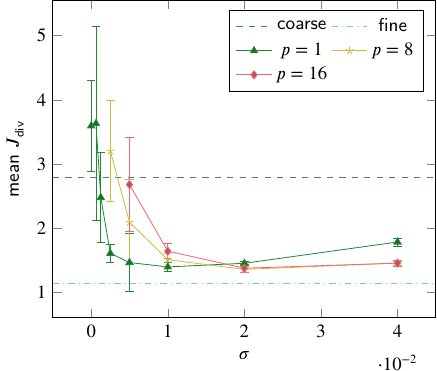}
	\fi
    \ifbuild
	\tikzsetnextfilename{noise_aug_ro1_vel}
	\begin{tikzpicture}
		\begin{axis}[%
			legend pos=north east,
			legend style={at={(0.97, 0.65)},anchor=east},
			xminorticks=true,
			xmin=-0.005,
			xmax=0.045,
			width=0.49\textwidth,
			xlabel={$\sigma$},
			ylabel={$e_{\bar{v}}$},
			]
			\addplot+[mark=none, dashed] table[x=setups, y=mean_v_err_coarse]{mean_gc1_noisy_data_reference.txt};
			\pgfplotsset{cycle list shift=1}
			
			\addplot+[error bars/.cd, y dir=both, y explicit] table[x=noise_sigma,y=mean, y error=std] {mean_gc1_noisy_data_pr1_mean_v_err_lines.txt};
			\addplot+[error bars/.cd, y dir=both, y explicit] table[x=noise_sigma,y=mean, y error=std] {mean_gc1_noisy_data_pr8_mean_v_err_lines.txt};
			\addplot+[error bars/.cd, y dir=both, y explicit] table[x=noise_sigma,y=mean, y error=std] {mean_gc1_noisy_data_pr16_mean_v_err_lines.txt};
			
			\legend{coarse, $p=1$, $p=8$, $p=16$};
		\end{axis}
	\end{tikzpicture}
    \else
	\includegraphics{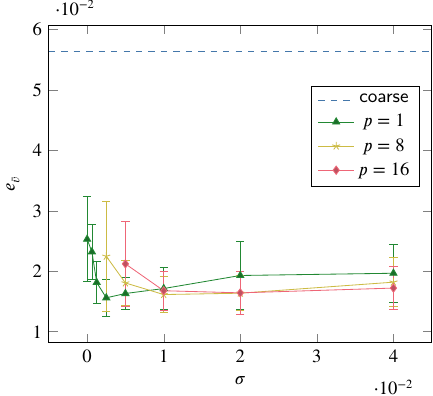}
	\fi
	\caption{The divergence $J_\text{div}$ \eqref{eq:div} and time-averaged velocity error $e_{\bar{v}}$ \eqref{eq:mean-vel-error} of NNs trained on noise augmented data. During training data generation the numerical solution is perturbed by adding noise sampled from a normal distribution with mean $\mu=0$ and standard deviation $\sigma$ to the velocity components with a probability of $c^{-1}$ per time-step ($c=1$ means every time-step). For reference, the mean magnitude of velocity components in the training set without perturbations is $\approx0.95$. Each point represents the mean from 25 NNs trained on the same dataset and the whiskers show the standard deviation. }
	\label{fig:div-noisy-aug}
\end{figure}
To determine the best parameters for $\sigma$ and $c$ we trained a number of small MLPs on datasets generated with different noise schedules. The training set consisted of the 3-round-obstacle case (ro3) while the evaluation was done on the single-obstacle case (ro1) with a different ellipse shape. For training, 1000 time-steps from ro3 and 600 time-steps each for validation from two variants of ro1 (with noise schedule $\sigma=0.04, c=1$ and without) where used for a total of \num{4.1e6} training samples and \num{2.4e6} validation samples.
% ro3: 4100096 samples, ro1: 1230848 (x2)

The results in \cref{fig:div-noisy-aug} clearly show the benefits of adding noise. Increasing the strength $\sigma$ significantly reduces the divergence of the NN solution and leads to much more consistent results. This effect is also visible in the time-averaged velocity, which is greatly improved by the NNs even with the instability but further improves with the added noise. For larger $\sigma$, the metrics get worse again as the training data is too far removed from the target task. A larger period $p$ dampens this effect and shifts the curve to the right. For $c=16$, the sweet spot appears to be around $\sigma=0.02$, or $2\%$ of the mean magnitude of the perturbed quantity. Consequently, we use a noise schedule with these parameters for subsequent experiments.

%%%%%%%%%%%%%%%%%%%%%%%%%%%%%%%%%%%%%%%%%%%%%%%%%%%%%%%%%%%%%%%%%%%%%%%%%%%%
\subsection{Rotation}
\begin{figure}
\centering
\ifbuild
\tikzsetnextfilename{rotation_ro1_freq}
\begin{tikzpicture}
    \pgfplotstableread{functionals_gc1_data_reference.txt}\loadedRefFunctionals
    \pgfplotstableread{functionals_gc1_rot_aug_freq_lines.txt}\loadedFunctionals
	\begin{axis}[
		symbolic x coords={ref,ref-rot90,small,small-rot90,large,large-rot90},
		% first one is skipped for some reason
		xticklabels={,ref,ref-\ang{90},small,small-\ang{90},large,large-\ang{90}},
		xmin={[normalized]-0.6},
		xmax={[normalized]+5.6},
		xticklabel style={
			tick label style={rotate=45},
		},
		xtick=data,
		ylabel={lift frequency $f$},%[\unit{s^{-1}}]
		legend pos=south east,
		width=0.495\textwidth,
		]
		
		\addplot+[mark=none, dashed] table[x=setups, y=freq_coarse]{\loadedRefFunctionals};
		\addplot+[mark=none, dashdotted] table[x=setups, y=freq_fine]{\loadedRefFunctionals};
		\foreach \x in {0,2,...,76} {
			\def\y{\the\numexpr\x+1\relax}
			\addplot table[x index=\x, y index=\y]{\loadedFunctionals};
		}
		\legend{coarse, fine};
	\end{axis}
\end{tikzpicture}
\else
\includegraphics{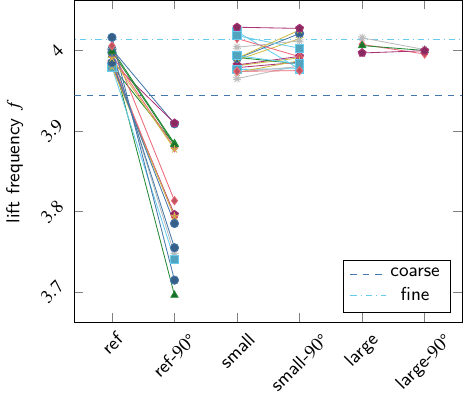}
\fi
\ifbuild
\tikzsetnextfilename{rotation_ro1_div}
\begin{tikzpicture}
    \pgfplotstableread{functionals_gc1_data_reference.txt}\loadedRefFunctionals
    \pgfplotstableread{functionals_gc1_rot_aug_mean_div_lines.txt}\loadedFunctionals
	\begin{axis}[
		symbolic x coords={ref,ref-rot90,small,small-rot90,large,large-rot90},
		% first one is skipped for some reason
		xticklabels={,ref,ref-\ang{90},small,small-\ang{90},large,large-\ang{90}},
		xmin={[normalized]-0.6},
		xmax={[normalized]+5.6},
		xticklabel style={
			tick label style={rotate=45},
		},
		ylabel={$J_\text{div}$},
		xtick=data,
		width=0.495\textwidth,
		]
	
		\addplot+[mark=none, dashed] table[x=setups, y=mean_div_coarse]{\loadedRefFunctionals};
		\addplot+[mark=none, dashdotted] table[x=setups, y=mean_div_fine]{\loadedRefFunctionals};
		\foreach \x in {0,2,...,76} {
			\def\y{\the\numexpr\x+1\relax}
			\addplot table[x index=\x, y index=\y]{\loadedFunctionals};
		}
		\legend{coarse, fine};
	\end{axis}
\end{tikzpicture}
\else
\includegraphics{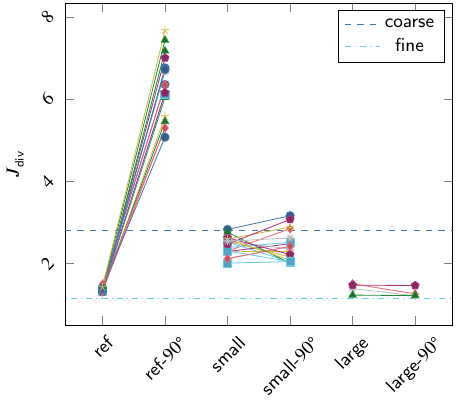}
\fi
    \caption{Functionals of the single obstacle case ro1 in the standard orientation and with the whole domain rotated by \ang{90} for simulations with different MLPs. Each line represents one MLP. The reference (ref), taken from \cref{sec:noise-aug}, is compared to MLPs trained on the same dataset but with rotation augmentation (small, large). Small networks have the same size as the reference (\num{1.3d+5} weights), large networks are roughly 8 times larger (\num{1.0d+6} weights).}
    \label{fig:rotaug}
\end{figure}

The computational problem we are trying to solve has rotational symmetry and the numerical solver is rotation-equivariant. The addition of the neural network in DNN-MG destroys this property. Although the local approach gives the neural network some ability to deal with different flow directions, simply rotating the domain by \qty{90}{\deg} in \cref{fig:rotaug} results in a new problem that proves too challenging when the neural network is not explicitly trained on similarly oriented flows. There is a large body of literature on imposing neural networks with equivariance; see~\cite{Rath2024} for a survey on existing methods. Here, we only consider data augmentation as a simple method, applicable to continuous symmetries,  to increase the robustness of DNN-MG.
More specifically, we extend the training procedure by rotating all vector-valued inputs and outputs (velocity components of $ \vec{x}_n, \vec{r}_n, \vec{d}_n $ and geometry $\vec{g}^*$) by a random angle $\alpha \in [0,2\pi)$ during batch creation. With augmented training, additional care is needed for the input and output normalization. The scaling and mean-shift also have to be rotation equivariant and so we group together components belonging to the same vector, set the shift to zero and assume zero-mean during the variance computation.

To test the effectiveness of rotation augmentation, we trained 25 MLPs each on the the 3-round-obstacles case (ro3) and evaluated them on the single-round-obstacle case (ro1) in the default orientation and rotated by \qty{90}{\deg}. In \cref{fig:rotaug} we see that the non-augmented MLPs (ref) disturb the simulation for the rotated case while the augmented MLPs (small, large) consistently improve the accuracy compared to the coarse solution. However, the augmented MLPs are far from equivariant as there is a significant dependence on the case and the solutions are overall worse than the baseline in the default orientation. The reduced skill is not a fundamental problem because it can be recovered by simply increasing the size of the MLP. We attribute this to the large MLP's greater capacity to memorize different orientations, as they still display a dependence on the orientation. Rotation augmentation is therefore a simple method to improve the robustness of the trained neural networks at the cost of increased training and inference time. Since the dominant flow in all our cases is in the same direction, they are unlikely to benefit from rotation augmentation and the following experiments are done without it. 

%%%%%%%%%%%%%%%%%%%%%%%%%%%%%%%%%%%%%%%%%%%%%%%%%%%%%%%%%%%%%%%%%%%%%%%%%%%%
\subsection{Multiple cases}\label{sec:train-data-dep}
While data augmentation allows us to get more out of an existing setup, we also have the option add more cases to generate more varied training data. To investigate the effect of different domains, we train 4 large MLPs each on noise augmented datasets derived from different cases. We compare training on only square obstacles (sq4, sq6), only round obstacles (ro3, ro5, ro9) and only small cases (sq4, ro3), with reference that includes all these cases. To ensure that each neural network receives the same amount of training, we effectively increase the number of training epochs for the smaller datasets. Two more cases (sq9, ro6) are held back for validation. For training on square obstacles only, just sq9 is used to select the best neural network, as the loss on ro6 tends to not improve at all.

\begin{table}
    \centering
    \begin{tabular}{lccccccccc}
\toprule
 & & \multicolumn{2}{c}{MSE loss sq9} & \multicolumn{2}{c}{MSE loss ro6} & \multicolumn{2}{c}{mean $\tau$~\eqref{eq:trunc-error} sq9} & \multicolumn{2}{c}{mean $\tau$~\eqref{eq:trunc-error} ro6} \\
 & train samples & mean & best & mean & best & mean & best & mean & best \\
\midrule
coarse & - & - & - & - & - & 0.045 & 0.045 & 0.023 & 0.023 \\
all & \num{3.1e7} & \textit{0.007} & \textit{0.007} & \textit{0.006} & \textit{0.006} & \textit{0.002} & \textit{0.002} & \textbf{0.002} & \textit{0.002} \\
sq4 ro3 & \num{7.8e6} & 0.088 & 0.086 & 0.016 & 0.015 & 0.007 & 0.007 & 0.004 & 0.003 \\
ro3 ro5 ro9 & \num{2.3e7} & 1.479 & 1.465 & \textbf{0.006} & \textbf{0.005} & 0.033 & 0.033 & \textit{0.002} & \textbf{0.002} \\
sq4 sq6 & \num{6.6e6} & \textbf{0.006} & \textbf{0.006} & 3.727 & 3.123 & \textbf{0.002} & \textbf{0.002} & 0.164 & 0.151 \\
\bottomrule
\end{tabular}
    \caption{Final validation loss and the local error of MLPs trained on different cases. The local error computed right after the fine warm-up at $n_0=50$ over $k=100$ steps for sq9 and $k=200$ steps for ro6}.
    \label{tab:train-dataset-dep}
\end{table}
Loss and a short-term local error on the two validation cases are recorded in \cref{tab:train-dataset-dep}. The evaluation time interval is short to circumvent stability issues that occur for some of the MLPs later on, which we address in \cref{sec:replay}. Unsurprisingly, more training cases are better and ``all'' performs well on both cases, on par with the neural networks trained only on one kind of domain. Furthermore, we can see that variety is key, as the small mixed dataset ``c4 gc3'' gives most of the improvement over the coarse simulation already. While the loss is $12$ times higher on sq9 when compared to the best neural networks, the local error is only factor $3.5$ higher and already tiny relative to the coarse baseline.
\begin{figure}
	% left bot right top
	\newcommand{\includeGraphicsCrop}[1]{\includegraphics[width=0.49\textwidth, trim={1.5cm 7cm 44cm 7cm}, clip]{#1}}
	\newlength{\labelDistance}
	\setlength{\labelDistance}{1.5cm}
    
    \centering
    \begin{subfigure}{0.492\textwidth}
    \ifbuild
    \tikzsetnextfilename{data_ro6}
    \begin{tikzpicture}[node distance=0.1cm, 
		label distance=-0.5cm,%
		]
        \tikzset{%
    	   pic/.style={inner sep=0pt},
        }
		\node(Coarse)[pic, label=north: coarse]{\includeGraphicsCrop{imgs/gc6/gc6-coarse.png}};
		\node(Square)[pic, right= of Coarse, label=north: sq4 sq6] {\includeGraphicsCrop{imgs/gc6/gc6-0_mlp_c4_c6.png}};
		\node(Fine)[pic, below= of Coarse, label=north: fine] {\includeGraphicsCrop{imgs/gc6/gc6-fine.png}};
		\node(Round)[pic, right= of Fine, label=north: ro3 ro5 ro9] {\includeGraphicsCrop{imgs/gc6/gc6-0_mlp_large_gc3_gc5_gc9.png}};

        \node(BoundingBox)[fit=(Coarse) (Round)] {};
		\node [pic, below= of BoundingBox]{\includegraphics[width=0.8\textwidth]{imgs/gc6/gc6-colorbar.png}};
	\end{tikzpicture}
    \else
    \includegraphics{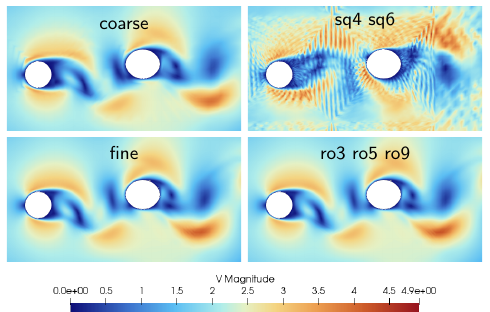}
    \fi
    \subcaption{Round obstacle (ro6)}
    \label{fig:train-dataset-dep-ro6}
    \end{subfigure}
    \renewcommand{\includeGraphicsCrop}[1]{\includegraphics[width=0.49\textwidth, trim={2.5cm 0.4cm 25cm 0.4cm}, clip]{#1}}
    \begin{subfigure}{0.492\textwidth}
    \ifbuild
    \tikzsetnextfilename{data_sq9}
    \begin{tikzpicture}[node distance=0.1cm, 
    	label distance=-0.45cm,%
    	]
        \tikzset{%
        	pic/.style={inner sep=0pt},
        }
    	\node(Coarse)[pic, label=north: coarse]{\includeGraphicsCrop{imgs/c9/c9-coarse.png}};
    	\node(Square)[pic, right= of Coarse, label=north: sq4 sq6] {\includeGraphicsCrop{imgs/c9/c9-c4_c6.png}};
    	\node(Fine)[pic, below= of Coarse, label=north: fine] {\includeGraphicsCrop{imgs/c9/c9-fine.png}};
    	\node(Round)[pic, right= of Fine, label=north: ro3 ro5 ro9] {\includeGraphicsCrop{imgs/c9/c9-gc3_gc5_gc9.png}};
    	
    	\node(BoundingBox)[fit=(Coarse) (Round)] {};
    	\node [pic, below= of BoundingBox]{\includegraphics[width=6cm]{imgs/c9/c9-colorbar.png}};
    \end{tikzpicture}
    \else
    \includegraphics{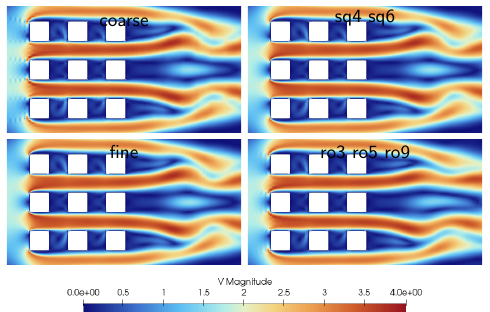}
    \fi
    \subcaption{Square obstacle (sq9)}
    \label{fig:train-dataset-dep-sq9}
    \end{subfigure}
    
    \caption{Velocity magnitude at $t=1.5$ ($n=150$) for two different cases. The neural networks used for DNN-MG where only trained on either square obstacles (sq4, sq6) or round obstacles (ro3, ro5, ro9). }
    \label{fig:train-dataset-dep}
\end{figure}
Trained just on square obstacles, DNN-MG ``sq4 sq6'' fails on the round obstacle case. The main problem is the shape of mesh cells. The square obstacle domains only consist of axis aligned rectangles and so the neural network can not deal with the trapezoid cells around round obstacles. This is easily visible in the velocity field, shown in \cref{fig:train-dataset-dep-ro6}, where artifacts appear tracing the grid. These artifacts already become visible after the first step with the neural network and can quickly destabilize the simulation. The other way around, training only on round obstacles, is more robust. There is still some improvement with ``ro3 ro5 ro9'' on sq9 over the coarse simulation, but the correction around the irregularities in front of the obstacles is muted in \cref{fig:train-dataset-dep-sq9}.

%%%%%%%%%%%%%%%%%%%%%%%%%%%%%%%%%%%%%%%%%%%%%%%%%%%%%%%%%%%%%%%%%%%%%%%%%%%%
\subsection{Replay Buffers}\label{sec:replay}
\begin{figure}
    \centering
    \ifbuild
    \tikzsetnextfilename{replay_ro9_local}
	\begin{tikzpicture}
        \pgfplotstableread{replay/V_gc9_ifix_mlp_replay_long_1step.txt}\loadedStepError
		\begin{axis}[
			xlabel={time $t [\unit{s}]$},
			ylabel={$\tau$},
			legend columns = 2,
			legend style={at={(0.55, 1.02)},anchor=south,font={\mystrut\footnotesize},
				legend cell align=left},
			every y tick scale label/.style={
				at={(0,1)},xshift=-10pt,yshift=5pt,anchor=south west,inner sep=0pt
			},
			width=0.495\textwidth,
			cycle multi list={
				bright\nextlist
				solid,dashed\nextlist
			},
			]
			% column 1 heading
			\addlegendimage{legend image with text=pred}
			\addlegendentry{}
			% column 2 heading
			\addlegendimage{legend image with text=interp}
			\addlegendentry{}
			
			\addlegendimage{legend image with text=}
			\addlegendentry{}
			\addplot+[mark repeat=64] table[x=time, y=coarse]{\loadedStepError};
			\addlegendentry{coarse}
			%skip dashed for coarse and fine color
			\pgfplotsset{cycle list shift=3};
			\foreach \nn/\nnlabel in {
				0_ifix_mlp_large_62/MLP,
				0_ifix_mlp_large_62_replay_15/MLP-replay%
%				1_ifix_mlp_large_58/MLP2,
%				1_ifix_mlp_large_58_replay_15/MLP2-replay%
			} {
				\addplot+[mark repeat=64] table[x=time, y=\nn, col sep=space]{\loadedStepError};
				\addlegendentry{}
				\addplot+[mark repeat=64] table[x=time, y=\nn-interp, col sep=space]{\loadedStepError};
				\addlegendentryexpanded{\nnlabel}
			}
			%\legend{coarse, MLP1, coarse-MLP1, MLP2, coarse-MLP2, MLP1-replay, coarse-MLP1-replay, MLP2-replay, coarse-MLP2-replay };
		\end{axis}
	\end{tikzpicture}
    \else
    \includegraphics{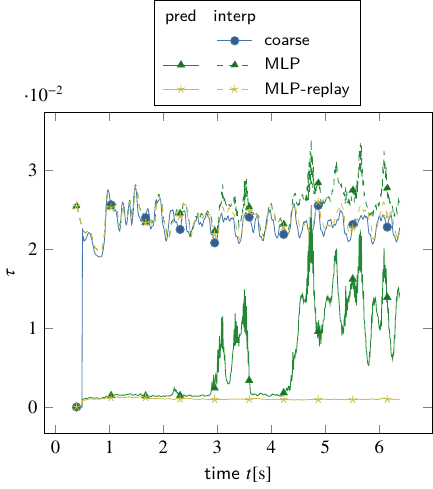}
    \fi
    \ifbuild
    \tikzsetnextfilename{replay_ro9_div}
	\begin{tikzpicture}
		\begin{axis}[
			xlabel={time $t [\unit{s}]$},
			ylabel={$J_{\text{div}}$},
			legend columns = 2,
			legend style={at={(0.5, 1.02)},anchor=south},
			width=0.495\textwidth,
			]
			\addplot+[mark repeat=64] table[x index=0, y index=1, col sep=space]{replay/fcts_gc9_ifix_mlp_replay_long_coarse_.txt};
			\addplot+[mark repeat=64] table[x index=0, y index=1, col sep=space]{replay/fcts_gc9_ifix_mlp_replay_long_fine.txt};
			\foreach \nn in {
				0_ifix_mlp_large_62,
				0_ifix_mlp_large_62_replay_15%
			} {
				\addplot+[mark repeat=64] table[x index=0, y index=1, col sep=space]{replay/fcts_gc9_ifix_mlp_replay_long_\nn_.txt};
			}
			\legend{coarse, fine, MLP, MLP-replay, MLP2, MLP2-replay };
		\end{axis}
	\end{tikzpicture}
    \else
    \includegraphics{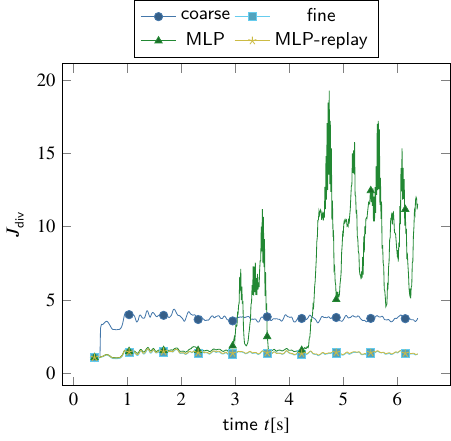}
    \fi
	\caption{Metrics of NN enhanced simulations compared with the reference simulations for ro9, the 9-round-obstacles case which is part of training set. The local error $\tau$~\eqref{eq:trunc-error} (left) is the $L_2$ difference between the velocity predicted by the NN and the numerical solution on the fine grid for the current time-step. Also provided is the error of the baseline, i.e. the coarse solution interpolated to the fine grid that the NN tries to correct. The cumulated effect on the simulation is visible in the divergence $J_{\text{div}}$~\eqref{eq:div} over time (right) with the standard reference using only the numerical solver. For a visualisation of the velocity with visible artifacts, see \cref{fig:artifacts-velocity}.}
	\label{fig:gc9-replay}
\end{figure}

Even with noise augmented training data, instabilities still develop in some cases with the DNN-MG simulation. In fact, the out-of-sample problem is most pronounced on the training case in \cref{fig:gc9-replay}, where the predictions are initially highly accurate, but produce massive deviations in the solution later on. A possible explanation is that the neural networks have a higher confidence on the training case, resulting in larger corrections. Since the NN is trained on minimizing the mean square error of a residual, the predictions will tend towards zero in uncertain situations.
When we consider only the residual for the current step, the NN prediction consistently improves upon the coarse solution, but when taken over many steps, the divergence gets much worse in the hybrid simulation. The fundamental problem is that there is a gap between the training and the inference task. Training samples are taken from a fine simulation with coarse inputs only computed for the current step. On the other hand, during inference, the coarse solver drives the simulation and we only approach the fine simulation if the neural network corrections are close to perfect.

To bridge the gap, we turn to replay buffers, a tool from reinforcement learning. We run the hybrid simulation with a trained NN in order to generate new trajectories for a second training stage. In each hybrid simulation step, the numerical solver is also used to compute a solution on the fine grid as ground truth, giving us a new training sample. Then the neural network is retrained on a combination of replay data, which teaches the NN how to correct its own errors, and the old data.

In theory, it should be best to generate small time-series and integrate the process of replay data generation tightly with the training since any updates to a neural network's weights will change the errors it produces. However, in \cref{fig:gc9-replay} we observe that significant deviations occur only after 250 time-steps of the full simulation. Thus, each replay iteration that adds valuable data likely comes at a considerable cost. Therefore, a process with two distinct stages consisting of a longer replay generation and a subsequent finetuning is much simpler to implement. In our experiments, this proved already effective for our method.

For this experiment, we first train large MLPs on a dataset consisting of multiple noise augmented simulations of square and round obstacle cases, simulated for 950 steps each. Other arrangements of obstacles with and without noise are retained as validation set. More details are given in \cref{appendix:experiment-details}. Afterward, each case from the training set is simulated again with the neural network corrections. To speed up the replay generation, which has to be done for every neural network, we shortened the simulation time from 950 steps to 650 steps and truncate the original data to balance the two datasets. Finally, both sets are combined to finetune the MLPs with the learning rate reduced from \num{0.01} to \num{0.0002}. Since the NNs are already trained, fewer epochs are needed to achieve convergence and we reduce the batch size from 512 to 256. %without risking instability during the training.
 
\begin{table}
    \centering
%	\scriptsize
	\begin{tabular}{lcccccccccc}
        \toprule
         & \multicolumn{2}{c}{mean $e_{J_d}$ \eqref{eq:drag-error}} & \multicolumn{2}{c}{mean $e_{J_l}$ \eqref{eq:lift-error}} & \multicolumn{2}{c}{mean $J_{\text{div}}$ \eqref{eq:div}} & \multicolumn{2}{c}{$e_{\bar{v}}$ \eqref{eq:mean-vel-error}} & \multicolumn{2}{c}{mean $\tau$ \eqref{eq:trunc-error}} \\
         & mean & best & mean & best & mean & best & mean & best & mean & best \\
        \midrule
        coarse & 0.044 & 0.044 & 0.087 & 0.087 & \textit{3.745} & \textit{3.745} & 0.050 & 0.050 & 0.023 & 0.023 \\
        baseline & \textit{0.008} & \textit{0.007} & \textit{0.013} & \textbf{0.010} & 11.076 & 9.894 & \textit{0.015} & \textit{0.012} & \textit{0.016} & \textit{0.014} \\
        replay & \textbf{0.007} & \textbf{0.005} & \textbf{0.010} & \textit{0.010} & \textbf{1.412} & \textbf{1.403} & \textbf{0.011} & \textbf{0.010} & \textbf{0.001} & \textbf{0.001} \\
        \midrule
        fine & 0.000 & 0.000 & 0.000 & 0.000 & 1.371 & 1.371 & 0.000 & 0.000 & 0.000 & 0.000 \\
        \bottomrule
    \end{tabular}
	\caption{Simulation results on ro9, a case from the training set, for 4 MLPs each. The drag $J_d$, lift $J_l$ and mean velocity error $e_{\bar{v}}$ are computed starting right after the fine warm-up at $n_0=50$ over $k=400$ steps. Divergence $J_{\text{div}}$ and local error $\tau$ are taken from $n_0=450$ over $k=500$ steps to show the longer-term behavior. The values are plotted over time in \cref{fig:replay-gc9}. }	
	\label{tab:replay-gc9}
\end{table}
\begin{table}
    \centering
%	\scriptsize
    \begin{tabular}{lcccccccccc}
    \toprule
     & \multicolumn{2}{c}{mean $J_d$ \eqref{eq:drag}} & \multicolumn{2}{c}{amp $J_l \eqref{eq:lift}$} & \multicolumn{2}{c}{freq \eqref{eq:freq}} & \multicolumn{2}{c}{mean $J_{\text{div}}$ \eqref{eq:div}} & \multicolumn{2}{c}{$e_{\bar{v}}$ \eqref{eq:mean-vel-error}} \\
     & mean & best & mean & best & mean & best & mean & best & mean & best \\
    \midrule
    coarse & 0.406 & 0.406 & 0.456 & 0.456 & 3.945 & 3.945 & 2.794 & 2.794 & 0.056 & 0.056 \\
    baseline & \textbf{0.486} & \textbf{0.484} & \textbf{0.576} & \textbf{0.568} & \textit{3.988} & \textit{3.995} & \textit{1.208} & \textit{1.200} & \textit{0.010} & \textit{0.009} \\
    replay & \textit{0.487} & \textit{0.486} & \textit{0.578} & \textit{0.575} & \textbf{3.989} & \textbf{3.996} & \textbf{1.163} & \textbf{1.153} & \textbf{0.008} & \textbf{0.006} \\
    \midrule
    fine & 0.474 & 0.474 & 0.558 & 0.558 & 4.013 & 4.013 & 1.144 & 1.144 & 0.000 & 0.000 \\
    \bottomrule
    \end{tabular}
\caption{Simulation results on ro1 for 4 MLPs each. All values are taken from around $n_0=490$ when the periodic flow is fully developed over $k=300$ steps. Each simulation is shifted slightly such that the lift period is in sync at the start. The lift amplitude is computed as $\text{amp} J_l = \max J_l - \min J_l$ over the interval. The values are plotted over time in \cref{fig:replay-gc1}.}
\label{tab:replay-gc1}
\end{table}
We present results for 4 random initialized MLPs. In \cref{tab:replay-gc9}, the replay training consistently prevents instabilities and also improves every metric, both in the short term and also for long simulations. The same is largely true for the validation case gc1, shown in \cref{tab:replay-gc1}, where the baseline is already stable. The drag and lift functionals show a bias as they further overshoot the reference. Interestingly, the functionals of the hybrid simulation get closer to the true solution, as determined by reference simulations on grids refined beyond the fine level. With two more levels we observe convergence towards amp $J_l=0.578$ and mean $J_d=0.493$. 
A possible explanation is that this is related to a difference between explicit and implicit methods, where, in cases of monotone convergence, they tend to approach the solution from opposite sides. The neural network with its purely additive update may act like an explicit method that counteracts the implicit solver's tendency to underestimate the functionals. 

\begin{figure}
    \centering
    \ifbuild
    \tikzsetnextfilename{replay_ro1_re}
	\begin{tikzpicture}
        \pgfplotstableread{replay/functionals_gc1_reynolds_mlp_mean_div_lines.txt}\loadedMeanDiv
		\begin{axis}[
			cycle list name=bright-grouped,
			xlabel={Reynolds number $\Rey$},
			ylabel={mean $J_{\text{div}}$},
			width=0.49\textwidth,
            legend style={at={(1.02, 0.5)},anchor=west},
			]
			% coarse
			\addplot+[mark repeat=1] table[x index=0, y index=17, col sep=space]{\loadedMeanDiv};
			% fine
			\addplot+[mark repeat=1] table[x index=0, y index=19, col sep=space]{\loadedMeanDiv};
			% baseline
			\foreach \y in {1,5,...,16} {
				\addplot+[mark repeat=1] table[x index=0, y index=\y, col sep=space]{\loadedMeanDiv};
			}
			\addlegendentry{baseline}
			% replay
			\foreach \y in {3,7,...,16} {
				\addplot+[mark repeat=1] table[x index=0, y index=\y, col sep=space]{\loadedMeanDiv};
			}
			\addlegendentry{replay}
			\legend{coarse, fine, baseline,,,, replay}
		\end{axis}
	\end{tikzpicture}
    \else
    \includegraphics{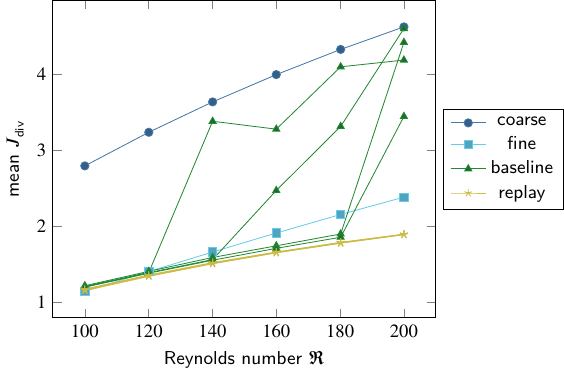}
    \fi
	\caption{Divergence $J_{\text{div}}$ on ro1 for 4 MLPs each, taken from $n_0=850$ over $k=200$ steps for simulations with increasing Reynolds number $\Rey$. The viscosity $\nu = 0.001$ is the same in all training cases, giving $\Rey=100$, and reduced to $\nu=0.0005$ for the highest Reynolds number.}
	\label{fig:replay-gc1-reynolds}
\end{figure}

Both ro1 and ro9 are close to the training data, but replay training also significantly improves the generalization ability of the NNs. All training cases have a Reynolds number around $\Rey = 100$. If we reduce the viscosity to increase $\Rey$, all the baseline MLPs quickly show increased divergence in \cref{fig:replay-gc1-reynolds} while the replay trained MLPs remain stable. The overall performance for higher $\Rey$ is investigated in \cref{sec:generalization}.

\begin{table}
    \centering
%    \scriptsize
	\begin{tabular}{lcccccc}
    \toprule
     & \multicolumn{2}{c}{mean $J_{\text{div}}$ \eqref{eq:div}} & \multicolumn{2}{c}{$e_{\bar{v}}$ \eqref{eq:mean-vel-error}} & \multicolumn{2}{c}{mean $\tau$ \eqref{eq:trunc-error}} \\
     & mean & best & mean & best & mean & best \\
    \midrule
    coarse & \textit{7.3196} & 7.3196 & 0.1248 & 0.1248 & 0.0452 & 0.0452 \\
    baseline & 9.9461 & \textit{5.2440} & \textit{0.0395} & \textbf{0.0282} & \textit{0.0173} & \textit{0.0036} \\
    replay & \textbf{5.1226} & \textbf{5.1082} & \textbf{0.0377} & \textit{0.0324} & \textbf{0.0029} & \textbf{0.0028} \\
    \midrule
    fine & 5.1242 & 5.1242 & 0.0000 & 0.0000 & 0.0000 & 0.0000 \\
    \bottomrule
\end{tabular}
	\caption{Simulation results on sq9* for 4 different MLPs. The 9-square-obstacle case sq9 is modified by adding uniform noise of strength $[-\num{5e-3}, \num{5e-3}]$, or up to \qty{5}{\%} of the square obstacle side length, to the interior vertices of the base mesh $\Omega_0$. All metrics are taken from $n_0=1050$ for $k=1000$ steps. The mean velocity does not fully converge in this case, see also \cref{fig:architectures-c9}. The much higher mean divergence of the baseline MLPs can largely be attributed to just one MLP, which introduces significant artifacts to the solution shown in \cref{fig:artifacts-velocity}.}
	\label{tab:replay-c9}
\end{table}
Another case that demonstrates the improved robustness is shown in \cref{tab:architectures-c9}. The square obstacle cases are difficult to evaluate because they lead to symmetric solutions that are unstable and highly sensitive. The symmetry is quickly broken by DNN-MG after 100-200 steps, and even inaccuracies in the implicit solver introduce deviations that can cause wildly different flows. We therefore randomly perturb the interior vertices of the base mesh, which results in a local geometry unseen during the training. Nevertheless, even the baseline MLPs perform well on this case, with only 1 out of the 4 displaying increased divergence. The replay training again addresses the remaining instability and leads to very consistent results. Beyond that, improvements are difficult to quantify because the flow has a complex vortex pattern and even averaged over 1000 steps there are still some fluctuations in the mean velocity field.

%%%%%%%%%%%%%%%%%%%%%%%%%%%%%%%%%%%%%%%%%%%%%%%%%%%%%%%%%%%%%%%%%%%%%%%%%%%%
\section{Generalization}\label{sec:generalization}
For DNN-MG to be effective, the cost of the training has to be amortized by the savings achieved with the hybrid simulations. It is therefore imperative that the neural networks generalize well, delivering improved accuracy on a wide range of problems. In this section, we evaluate the impact of the neural network architecture and the value of additional information provided to the neural network. 

We train 4 large neural networks each of the MLP, RNN and Transformer architectures and another set of MLPs with patch size $M=1$ instead of $0$. Additionally, we consider 2 sets of MLPs which operate on a coarser level with $J=2$ instead of $J=1$. As these start with a much worse solution, they are not competitive and we only give broad results in the following and refer to \cref{appendix:jmplvl} for the details. The number of learnable weights is around \num{1e6} for all networks except for the Transformers, which have only \num{6e5} weights but require significantly more compute per weight. The training procedure is the same as described in \cref{sec:replay}. The training data, a mixture of round and square obstacle cases, is generated with added noise. Then each NN is fine-tuned on the training set, augmented with a single generation of replays. No special treatment is needed for the Transformers to stabilize this second training and the same batch-size and constant learning rate are used for all architectures.

%%%%%%%%%%%%%%%%%%%%%%%%%%%%%%%%%%%%%%%%%%%%%%%%%%%%%%%%%%%%%%%%%
\subsection{Accuracy}
\begin{figure}
	\begin{subfigure}{0.495\textwidth}
    \ifbuild
    \tikzsetnextfilename{accuracy_loss_ro6}
	\begin{tikzpicture}
        \pgfplotstableread[col sep=comma]{loss_valid_gc6.txt}\loadedData
		\begin{axis}[
			cycle list name=bright-grouped,
			xlabel={epoch},
			ylabel={MSE loss},
			width=0.99\textwidth,
			legend style={at={(0.97, 0.45)},anchor=east},
			]
			\pgfplotsset{cycle list shift=2}
			% MLP
			\foreach \y in {2,5,...,12} {
				\addplot+[mark repeat=1] table[x index=0, y index=\y, col sep=comma]{\loadedData};
			}
            \pgfplotsset{cycle list shift=6}
            % GRU
			\foreach \y in {1,4,...,12} {
				\addplot+[mark repeat=1] table[x index=0, y index=\y, col sep=comma]{\loadedData};
			}
			% TF
            \foreach \y in {3,6,...,12} {
				\addplot+[mark repeat=1] table[x index=0, y index=\y, col sep=comma]{\loadedData};
			}
			\legend{MLP,,,,RNN,,,,Transformer}
		\end{axis}
	\end{tikzpicture}
    \else
    \includegraphics{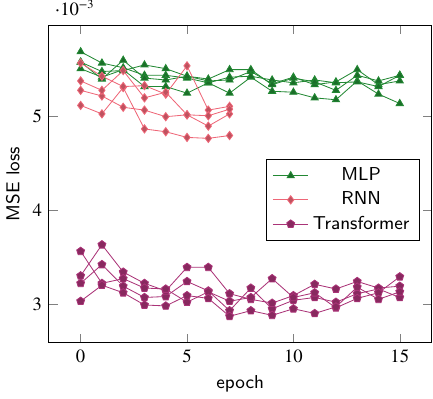}
    \fi
	\subcaption{Round obstacle (ro6)}
	\end{subfigure}
	\begin{subfigure}{0.495\textwidth}
        \ifbuild
        \tikzsetnextfilename{accuracy_loss_c9}
		\begin{tikzpicture}
            \pgfplotstableread[col sep=comma]{loss_valid_c9.txt}\loadedData
			\begin{axis}[
				cycle list name=bright-grouped,
				xlabel={epoch},
				ylabel={MSE loss},
				width=0.99\textwidth,
				]
				\pgfplotsset{cycle list shift=2}
				% MLP
				\foreach \y in {2,5,...,12} {
					\addplot+[mark repeat=1] table[x index=0, y index=\y, col sep=comma]{\loadedData};
				}
                \pgfplotsset{cycle list shift=6}
				% GRU
				\foreach \y in {1,4,...,12} {
					\addplot+[mark repeat=1] table[x index=0, y index=\y, col sep=comma]{\loadedData};
				}
                % TF
				\foreach \y in {3,6,...,12} {
					\addplot+[mark repeat=1] table[x index=0, y index=\y, col sep=comma]{\loadedData};
				}
				\legend{MLP,,,,RNN,,,,Transformer}
			\end{axis}
		\end{tikzpicture}
        \else
        \includegraphics{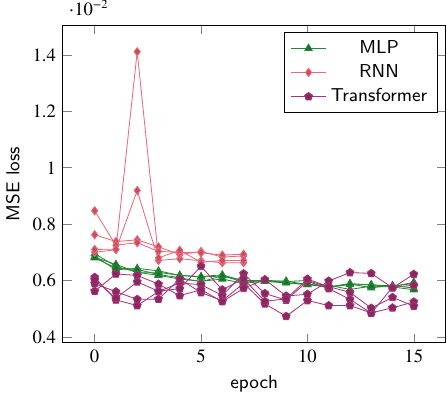}
        \fi
		\subcaption{Square obstacles (sq9)}
	\end{subfigure}
	\caption{Validation loss of different cases during replay training for 4 NNs of each architecture.}
	\label{fig:valid-loss-architectures}
\end{figure}
Already during training, some visible differences emerge between the architectures. The validation loss, visualized in \cref{fig:valid-loss-architectures}, shows a clear gap between the Transformers and the other architectures on round obstacle cases. The final loss after convergence is consistently \qty{40}{\%} lower for the Transformers on ro6. On the other hand, on the square obstacle case sq9 the differences are less pronounced and the RNNs seem to perform slightly worse. However, the validation loss is not necessarily an accurate indicator of later performance in DNN-MG due to the distributional shift explored in \cref{sec:noise-aug}. It is also not possible to directly compare values with different patch sizes because the occurrence of repeated nodes at the patch boundaries changes.

\begin{table}
    \centering
	\begin{tabular}{lcccccccc}
    \toprule
     & \multicolumn{2}{c}{mean $e_{J_d}$ \eqref{eq:drag-error}} & \multicolumn{2}{c}{mean $e_{J_l}$ \eqref{eq:lift-error}} & \multicolumn{2}{c}{mean $J_{\text{div}}$ \eqref{eq:div}} & \multicolumn{2}{c}{mean $\tau$ \eqref{eq:trunc-error}} \\
     & mean & best & mean & best & mean & best & mean & best \\
    \midrule
    coarse & 0.0308 & 0.0308 & 0.0976 & 0.0976 & 3.3068 & 3.3068 & 0.0235 & 0.0235 \\
    MLP & \textit{0.0077} & \textit{0.0061} & \textit{0.0251} & \textit{0.0189} & 1.3476 & 1.3374 & 0.0019 & 0.0018 \\
    MLP-M1 & 0.0081 & 0.0066 & 0.0267 & 0.0211 & \textbf{1.3146} & \textit{1.3115} & \textbf{0.0014} & \textbf{0.0014} \\
    RNN & 0.0111 & 0.0067 & 0.0367 & 0.0227 & 1.3495 & 1.3421 & 0.0019 & 0.0019 \\
    Transf & \textbf{0.0068} & \textbf{0.0055} & \textbf{0.0214} & \textbf{0.0152} & \textit{1.3171} & \textbf{1.3106} & \textit{0.0015} & \textit{0.0015} \\
    \midrule
    fine & 0.0000 & 0.0000 & 0.0000 & 0.0000 & 1.2527 & 1.2527 & 0.0000 & 0.0000 \\
    \bottomrule
    \end{tabular}
	\caption{Simulation results on ro6 for 4 NNs each of different architectures. All metrics are computed starting right after the warmup at $n_0=50$. Drag error $e_{J_d}$ and lift error $e_{J_l}$ are taken over $k=200$ steps to keep the influence of the temporal shift small. Divergence $J_{\text{div}}$ and the local error $\tau$ are averaged over $k=600$ steps. The values are plotted over time in \cref{fig:architectures-gc6}.}
	\label{tab:architectures-gc6}
\end{table}
A look at actual simulation results in \cref{tab:architectures-gc6} shows the same ranking (Transformer $>$ MLP $>$ RNN) as the validation loss but the differences are much smaller. All neural networks show a significant improvement over the coarse simulation with metrics improving by factors ranging from $2.5$ to $16.7$. The increased patch size of MLP-M1 makes them compete with Transformers for the global metrics, but the functionals indicate that the correction close to the obstacle is a little worse compared to patch size $M=0$.

\begin{table}
    \centering
%	\scriptsize
	\begin{tabular}{lcccccccc}
    \toprule
     & \multicolumn{2}{c}{$\min J_l$ \eqref{eq:lift}} & \multicolumn{2}{c}{$\max J_l$ \eqref{eq:lift}} & \multicolumn{2}{c}{mean $J_{\text{div}}$ \eqref{eq:div}} & \multicolumn{2}{c}{$e_{\bar{v}}$ \eqref{eq:mean-vel-error}} \\
     & mean & best & mean & best & mean & best & mean & best \\
    \midrule
    coarse & -0.0263 & -0.0263 & -0.0001 & -0.0001 & 7.3186 & 7.3186 & 0.1248 & 0.1248 \\
    MLP & \textit{-0.0793} & \textbf{-0.0752} & \textbf{0.0511} & \textbf{0.0516} & \textit{5.1221} & \textit{5.1081} & 0.0377 & 0.0324 \\
    MLP-M1 & -0.0936 & \textit{-0.0756} & 0.0838 & 0.0524 & 8.5799 & 5.2860 & 0.0492 & \textit{0.0238} \\
    RNN & \textbf{-0.0779} & -0.0743 & \textit{0.0530} & 0.0520 & \textbf{5.0278} & \textbf{4.9998} & \textit{0.0300} & 0.0246 \\
    Transf & -0.0821 & -0.0783 & 0.0486 & \textit{0.0512} & 5.2428 & 5.2184 & \textbf{0.0246} & \textbf{0.0206} \\
    \midrule
    fine & -0.0750 & -0.0750 & 0.0514 & 0.0514 & 5.1255 & 5.1255 & 0.0000 & 0.0000 \\
    \bottomrule
    \end{tabular}
	\caption{Simulation results on sq9* for 4 NNs each of different architectures. The 9-square-obstacle case sq9 is modified by adding uniform noise of strength $[-\num{5e-3}, \num{5e-3}]$, or up to \qty{5}{\%} of the square obstacle side length, to the interior vertices of the base mesh $\Omega_0$. All metrics are taken from $n_0=1050$ for $k=1000$ steps. The MLP results are the same as ``replay'' in \cref{tab:replay-c9}. The values are plotted over time in \cref{fig:architectures-c9}.}
	\label{tab:architectures-c9}
\end{table}
To evaluate the performance on sharp corners, we again turn to the case sq9* with a randomly perturbed mesh from \cref{sec:replay}. This is a more challenging test because there are no such examples of modified geometry in the training set. Again, the neural networks deliver significant improvements over the coarse solution. The amplification of the lift $J_l$ is of note because, while still small in absolute terms, it is a significant deviation from the square obstacles in the training data where $J_l \approx 0.0$. Between the different architectures there is no clear winner on sq9* as the flow is too chaotic to measure their minute differences. Only for MLP-M1 2 out of the 4 neural networks display significant deviations that are visible in the mean velocity over time in \cref{fig:architectures-c9}, and one MLP-M1 has increased divergence, a sign of instability.

\begin{figure}
\centering
\ifbuild
\else
\includegraphics{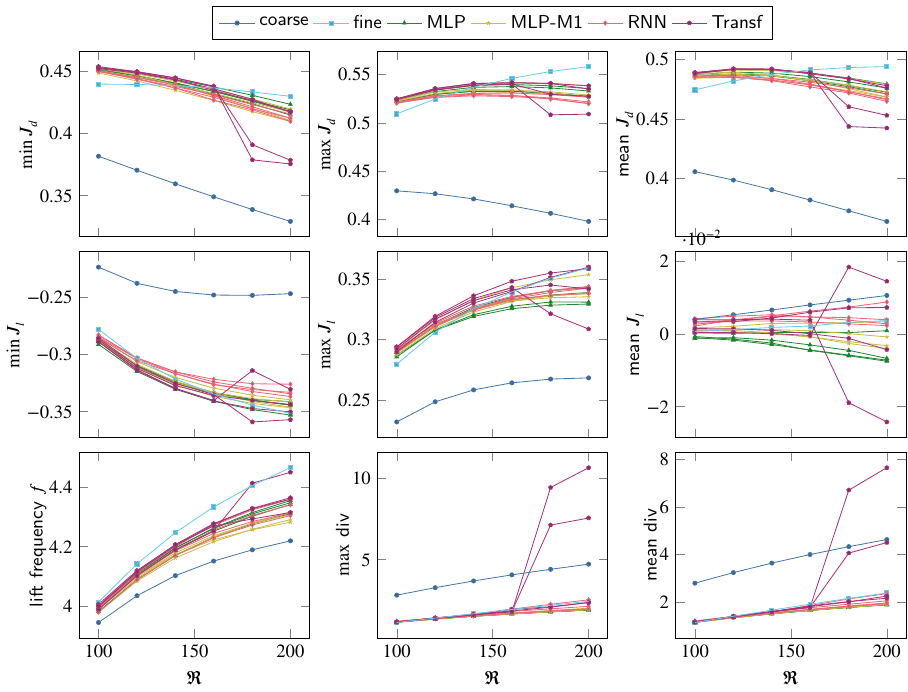}
\fi
	\caption{Functionals on ro1 for 4 NNs each, taken from $n_0=850$ over $k=200$ steps for simulations with increasing Reynolds number $\Rey$. The viscosity $\nu = 0.001$ is the same in all training cases, giving $\Rey=100$, and reduced to $\nu=0.0005$ for the highest Reynolds number.}
	\label{fig:gc1-architectures}
\end{figure}
In \cref{sec:replay} we observed that replay training makes the neural networks more robust to changes in the Reynolds number. The same setup on ro1 is extended to the architecture comparison in \cref{fig:gc1-architectures}. Broadly, the simulation remains stable even when doubling $\Rey$ and the neural networks provide corrections of a similar magnitude as the fine reference simulation. For $\Rey=100$ the functionals are systematically over-estimated, with all NNs displaying the same bias as the MLPs in \cref{tab:replay-gc1}. Again, the values are closer to those of a reference simulation of a higher level, though the lift overshoots in some NN simulations even compared to that reference. The largest corrections in terms of magnitude come from the Transformers, followed by the MLPs and they tend to be more accurate at higher $\Rey$. However, the stronger predictions also increase the risk of NN instability, as 2 out of the 4 Transformers clearly show in every functional.
This aside, the most interesting functionals are perhaps $\min J_l$ and $\max J_d$, since we see qualitatively different trends in the coarse and fine simulation. While $J_d$ decreases in the coarse solution as $\Rey$ becomes larger, it increases in the fine simulation and $J_l$ levels out in the coarse simulation but keeps decreasing in the fine simulation. Both trends are reproduced by the NNs to some extent. 

\begin{table}
    \centering
%	\scriptsize
	\begin{tabular}{lcccccccc}
        \toprule
         & \multicolumn{2}{c}{mean $J_d$ \eqref{eq:drag}} & \multicolumn{2}{c}{amp $J_d$ \eqref{eq:drag}} & \multicolumn{2}{c}{mean $J_{\text{div}}$ \eqref{eq:div}} & \multicolumn{2}{c}{$e_{\bar{v}}$ \eqref{eq:mean-vel-error}} \\
         & mean & best & mean & best & mean & best & mean & best \\
        \midrule
        coarse & 0.3484 & 0.3484 & 0.0487 & 0.0487 & 4.7263 & 4.7263 & 0.0516 & 0.0516 \\
        MLP & 0.4126 & 0.4116 & \textit{0.0656} & \textit{0.0651} & \textbf{3.1876} & \textbf{3.0969} & 0.0249 & 0.0220 \\
        MLP-M1 & \textit{0.4121} & \textit{0.4111} & 0.0660 & 0.0655 & 3.4354 & 3.2912 & \textbf{0.0231} & \textit{0.0216} \\
        RNN & 0.4124 & 0.4112 & 0.0668 & 0.0657 & \textit{3.2376} & \textit{3.1435} & 0.0246 & 0.0225 \\
        Transf & \textbf{0.4088} & \textbf{0.4014} & \textbf{0.0652} & \textbf{0.0650} & 3.5942 & 3.4090 & \textit{0.0238} & \textbf{0.0197} \\
        \midrule
        fine & 0.3986 & 0.3986 & 0.0625 & 0.0625 & 2.9192 & 2.9192 & 0.0000 & 0.0000 \\
        \bottomrule
    \end{tabular}
	\caption{Simulation results on ro3* for 4 NNs each of different architectures. Compared to ro3, the radii of the ellipses are changed and the middle of the three successive obstacles is kept square during mesh refinement. All metrics are computed starting at $n_0=400$ and computed over $k=600$ steps. The amplitude of the drag is computed as $\text{amp} J_d = \max J_d - \min J_d$. The values are plotted over time in \cref{fig:architectures-sgc3}.}
	\label{tab:architectures-sgc3}
\end{table}
Finally, we consider a case ro3* that combines both round and square obstacles. As base serves ro3, three round obstacles arranged one after the other in the channel. Although the original case is part of the training set, we change the radii of the first and third obstacle while leaving the middle one a square during the mesh refinement. This leads to a significant change in the flow behavior. The changes are both local in that we have a square obstacle with the mesh topology of a round obstacle and non-local in the interaction between round and square obstacles. The simulation results are summarized in \cref{tab:architectures-sgc3}. Once more, all neural networks produce excellent results but Transf and MLP-M1 take the lead overall. Again, a more accurate mean velocity $e_{\bar{v}}$ correlates with increased divergence, but all DNN-MG solutions are well below the divergence of the coarse solution.

The $J=2$ variants of the MLPs, recorded in \cref{appendix:jmplvl}, perform strictly worse than the $J=1$ variants across all experiments. Nevertheless, the hybrid simulations running with a coarser numerical solver are stable and the local-in-time metrics $\tau$ and $J_\text{div}$ are improved well beyond one level of the numerical solver. For long term metrics and the functionals the picture is mixed with overall results closer to the numerical solution on the intermediate level, i.e. the coarse solution for the $J=1$ variants. A challenge is that the coarser grid does not have enough degrees of freedom for the typical dynamics to even develop in the numerical solution. Also, given that the reconstruction task is more difficult with $J=2$, ideally a larger neural network should be used to achieve satisfactory results.

%%%%%%%%%%%%%%%%%%%%%%%%%%%%%%%%%%%%%%%%%%%%%%%%%%%%%%%%%%%%%%%%%
\subsection{Performance}
Two aspects determine the performance of DNN-MG, the initial cost of the training and the later savings achievable during hybrid simulations. For the latter, a theoretical justification for DNN-MG was already given in~\cite{MargenbergRichter2020,Margenberg2024} based on the reduction in required floating point operations. In practice, an important advantage of DNN-MG is that it enables the effective utilization of GPUs. Of course, one can also port the numerical solver itself to GPU as was done for Gascoigne 3D, the library also used in this work, in~\cite{Liebchen2024}.
However, the implementation can be challenging and not every algorithm is suitable for a GPU implementation. Sparse solvers generally achieve a low device utilization as they are heavily bandwidth limited. Deep neural networks, on the other hand, can take full advantage of the available compute, both because of the algorithms used and because modern GPUs are designed for deep learning. Since the number of unknowns grows geometrically across multigrid levels, the two finest levels handled by the neural network account for the vast majority of all degrees of freedom — 75\% in 2D and over 90\% in 3D — making the neural network correction the computationally dominant component of DNN-MG. For two jump levels ($J=2$), more than $90\%$ of the unknowns in 2D and about $98\%$ in 3D are neural network contributions.

Since the focus of this work is not on performance, both simulation speed and training were only tuned to an acceptable level to conduct the experiments. Further optimizations are likely possible. For example, the Newton-Krylov geometric multigrid method we use may not be the most efficient in practice~\cite{ahmedAssessmentSolversSaddle2018}, but DNN-MG could be used even with single level solvers. 

\begin{table}
    \centering
%    \scriptsize
	\begin{tabular}{cccccc}
		\toprule
		& MLP [h] & MLP-M1 [h] & MLP-J2 [h]  & RNN  [h] & Transf [h] \\
		\midrule
        data generation & 5.5 & 5.2 & 4.9 & 5.5 & 5.5 \\
		initial training & 8.4 &  $9.4^*$ & $11.6^*$ & 17.4 & 12.3 \\
        replay generation & 2.1 & 2.0 & 2.0 & 2.2 & 2.1 \\
		replay training & 9.5 & $6.0^*$ & $6.5^*$ & 16.0 & 19.8 \\
        \midrule
        total & 25.5 & 22.6 & 25.0 & 41.1 & 39.7 \\
		\bottomrule
	\end{tabular}
	\caption{Training times on a H100 PCIe GPU for large neural networks with \num{1e6} weights (MLP, MLP-M1, MLP-J2, RNN) and \num{6e5} weights (Transformer), respectively. Initial data generation, replay generation were done on a system with a RTX 3090 GPU and a i9-10900X CPU and the variations in timings between the architectures are mostly background noise; see \cref{fig:runtime} for comparison. For the initial training, the data can be reused except when the patch size ($M$) or jump level ($J$) changes. Timings for replay generation and training are per NN. Training runs marked with a ``$^*$'' are from the RTX 3090 system and feature a smaller validation set, with purely diagnostic cases removed. The differences between the MLPs can be largely attributed to these changes and background noise.}
	\label{tab:train-time}
\end{table}
The cost of training is recorded in \cref{tab:train-time} for each architecture. MLP and Transformers were trained for the same number of epochs, MLP-M1 for 4 times as many to get the same total number of weight updates as the base MLP. For the RNN, the number of epochs had to be reduced because backpropagation through time with long roll-outs is much more expensive. The simpler MLP architecture demonstrates a clear advantage here, as training was twice as fast as those for the more complicated architectures. However, for every architecture,
%for every architecture training was well within the regime of diminishing returns and 
both data efficiency and training time could likely be significantly improved without much loss in accuracy. While the validation loss continues to improve with more epochs, the gains in simulation skill are probably minor. In small scale tests with checkpoints of MLPs taken after a little more than half the training epochs, we barely see a difference in the hybrid simulation. Furthermore, during validation the neural networks were evaluated on multiple cases purely as diagnostics. As these cases did not factor into the selection of the best neural network, their evaluation could be skipped entirely without affecting the outcome, reducing the total number of samples used for training and validation by up to $1/3$.

\begin{figure}
	\begin{subfigure}{0.325\textwidth}
    \ifbuild
    \tikzsetnextfilename{performance_ro1}
	\begin{tikzpicture}
        \pgfplotstableread{goodcheese1_timings.txt}\loadedData
		\begin{axis}[
			ybar stacked,
            label style={font=\scriptsize},
		      tick label style={font=\scriptsize},
			ymin=0.0,%
			ymajorgrids=true,
            ylabel={total time $[s]$},
			symbolic x coords={coarse, fine, MLP, MLP-M1, MLP-J2, RNN, Transf},
			xtick=data,
            tick label style={rotate=35},
            legend style={legend pos=north east,font=\footnotesize},
			width=0.99\textwidth,
            bar width=0.3cm,
			]
			\addplot[stack plots=false, error bars/.cd,y dir=both,y explicit] table [x=file, y=sum-mean, y error=sum-std]{\loadedData};
			\pgfplotsset{cycle list shift=-1}
			\addplot table [x=file, y=Solver-mean, y error] {\loadedData};
			\addplot table [x=file, y=other-mean] {\loadedData};
			\addplot table [x=file, y=Network-mean] {\loadedData};
			\legend{,Solver, Interface, Neural network};
		\end{axis}
	\end{tikzpicture}
    \else
    \includegraphics{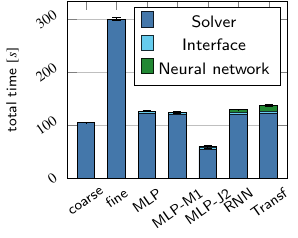}
    \fi
	\subcaption{ro1}
	\label{fig:runtime-gc1}
	\end{subfigure}
	\begin{subfigure}{0.325\textwidth}
        \ifbuild
        \tikzsetnextfilename{performance_c9}
		\begin{tikzpicture}
            \pgfplotstableread{cheese9_5e-3r_timings.txt}\loadedData
			\begin{axis}[
				ybar stacked,
                label style={font=\scriptsize},
		          tick label style={font=\scriptsize},
                enlarge y limits=0.15,
				ymin=0.0,%
                ymajorgrids=true,
                ylabel={total time $[s]$},
				symbolic x coords={coarse, fine, MLP, MLP-M1, MLP-J2, RNN, Transf},
                tick label style={rotate=35},
				xtick=data,
				legend style={legend pos=north east,font=\footnotesize},
				width=0.99\textwidth,
                bar width=0.3cm,
				]
				\addplot[stack plots=false, error bars/.cd,y dir=both,y explicit] table [x=file, y=sum-mean, y error=sum-std]{\loadedData};
				\pgfplotsset{cycle list shift=-1}
				\addplot table [x=file, y=Solver-mean] {\loadedData};
				\addplot table [x=file, y=other-mean] {\loadedData};
				\addplot table [x=file, y=Network-mean] {\loadedData};
				\legend{,Solver, Interface, Neural network};
			\end{axis}
		\end{tikzpicture}
        \else
        \includegraphics{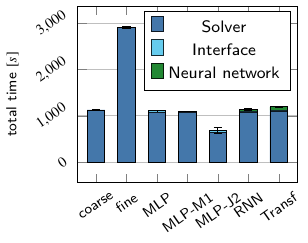}
        \fi
		\subcaption{sq9*}
		\label{fig:runtime-c9r}
	\end{subfigure}
    \begin{subfigure}{0.325\textwidth}
        \ifbuild
        \tikzsetnextfilename{performance_its}
        \begin{tikzpicture}
        \pgfplotstableread{newton_its.txt}\loadedData
    	\begin{axis}[
            ybar,
            tick align=inside,
            label style={font=\scriptsize},
		      tick label style={font=\scriptsize},
            ymin=0.0,%
            ymajorgrids=true,
            ylabel={average \#iterations},
            symbolic x coords={coarse, fine, MLP, MLP-M1, MLP-J2, RNN, Transf},
            tick label style={rotate=35},
            xtick=data,
            legend style={legend pos=south east,font=\small},
            width=0.99\textwidth,
            bar width=0.18cm,
            ]
            %\addplot[stack plots=false, error bars/.cd,y dir=both,y explicit] table [x=file, y=sum-mean, y error=sum-std]{\loadedData};
            \addplot+[error bars/.cd,y dir=both] table [x=name, y=ro1-mean, y error=ro1-std] {\loadedData};
            \addplot+[error bars/.cd,y dir=both] table [x=name, y=sq9s-mean, y error=Sq9s-std] {\loadedData};
            \legend{ro1, sq9*};
    	\end{axis}
        \end{tikzpicture}
        \else
        \includegraphics{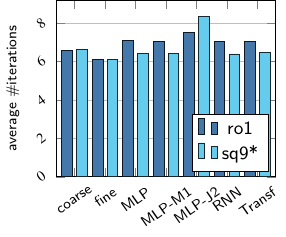}
        \fi
        \subcaption{Newton iterations}
        \label{fig:runtime-newtonits}
    \end{subfigure}
	\caption{Runtime of the simulation for different cases. The bars show the mean from 4 runs and the whiskers indicate the standard deviation. The warmup of 50 steps with just the numerical solver is not measured since this work is the same in every run. Neural network inference is done on a RTX 3090 GPU while the numerical solver runs with 10 threads on a i9-10900X CPU. In ``Interface'', additional work done by the solver which is needed to interface with the NN is collected. This includes fine residual computation, gathering of local inputs and the addition of the outputs to the global solution. Everything that happens in PyTorch, i.e. host-device memory transfers, normalization and additional input preparation for RNNs and Transformers are attributed to the category ``Neural Network''. }
	\label{fig:runtime}
\end{figure}
More important than training is the simulation runtime of DNN-MG, which we present in \cref{fig:runtime} for two different cases. In general, the numerical solver dominates the runtime. The inference time for the MLPs is negligible, but even for the 6 times slower RNNs and the 12 times slower Transformers neural network inference is not a major factor in the runtime of the simulation. On ro1, the NN enhanced simulation is thereby around $2.2$ times faster than the fine solution. An increase of \qty{20}{\%}-\qty{30}{\%} compared to the coarse simulation can be attributed to additional 0.5 Newton steps, shown in \cref{fig:runtime-newtonits}, that are needed on average to reach the same tolerance for the NN perturbed solution. On sq9*, the opposite effect can be seen, where on average 0.2 fewer Newton steps are needed by the hybrid simulations, causing MLP and MLP-M1 in particular to be slightly faster than just the coarse simulation. The speedup over the fine simulation is around $2.6$.

Since the speedup is mainly determined by the difference between the coarse and fine solver, much greater gains are possible for more difficult problems and with more mesh refinement levels in between coarse and fine grids. Starting one level lower and going to $J=2$ (cf. \cref{fig:patches}), MLP-J2 is twice as fast as $J=1$ simulations in \cref{fig:runtime}. The speedup of the solver is less for this level because the coarse solution is too small, having just 1728 (ro1) and 4368 (sq9*) degrees of freedom, respectively. For the 3D Navier-Stokes equations we achieved a speedup of factor 35 with a jump level $J=2$~\cite{Margenberg2024}. Thus, the main takeaway from the present results is in the difference in scaling between the numerical solver and the neural networks. The fine solutions have, respectively, 25344 (ro1) and 64776 (sq9*) degrees of freedom, while the neural networks have on the order of \num{1e6} weights. This demonstrates that increasing the size of the neural network to achieve higher accuracy is possible without a significant impact on the overall runtime of the method, in contrast to regular mesh refinement in the numerical solver, which incurs a considerable cost.

Given the cost of inference is not a major factor during the simulation, even larger MLPs would be feasible for the given setup. However, for Transformers and RNNs the cost could be more of a limiting factor. It is therefore worth noting that the neural networks in our implementation use the raw C++ PyTorch API which prevents optimizations across operations that are commonly done by other inference solutions like TorchInductor. Transformers in particular are the focus of current hardware and software development, e.g. \cite{TransformerEnigneGit}, since they are ubiquitous in Generative AI~\cite{Sengar2024}. Thus, the neural network component of DNN-MG could be made more computationally efficient if needed.

%%%%%%%%%%%%%%%%%%%%%%%%%%%%%%%%%%%%%%%%%%%%%%%%%%%%%%%%%%%%%%%%%
\subsection{Limitations}
DNN-MG is ill-suited for problems with highly sensitive solutions. This limitation was already discussed in the context of the square obstacle cases, which have a symmetric solution that constitutes an unstable equilibrium. The neural networks always introduce a small but unpredictable error that, if amplified, can produce different macroscopic flows. In contrast to a conventional numerical scheme, it is currently not possible to validate a solution by increasing the mesh resolution and using a smaller time-step.

Another problem we encountered in rare occasions during testing was an interaction between the do-nothing outflow boundary and the NN corrections. A spurious inflow generated by the NN can be amplified over multiple steps, thus destabilizing the whole simulation. This issue is distinct from the out-of-sample behavior or NN instability, since the error is highly localized and not necessarily accompanied by a visible increase in the divergence. This problem is well known and one possible solution is to use the directional-do-nothing condition~\cite{BraackMucha2014}, which stabilizes backflow into the domain.  However, additional (replay) training also helps to minimize this problem when the classical do-nothing condition is used.

%%%%%%%%%%%%%%%%%%%%%%%%%%%%%%%%%%%%%%%%%%%%%%%%%%%%%%%%%%%%%%%%%%%%%%%%%%%%
\section{Conclusion}\label{sec:conclusion}
We have presented improvements to the deep neural network multigrid solver (DNN-MG), a hybrid finite element / deep learning method for the efficient solution of boundary value problems. Using the example of the nonstationary 2D Navier-Stokes equation, we demonstrated improved accuracy and generalizability over our previous work by optimizing the neural networks and extending the training procedure with data augmentation.

We found that an increased spatial receptive field helps the neural network in the solution of the non-local Navier-Stokes equations, while an extension to multiple past time-steps has little effect. To make use of multiple cells on an unstructured mesh without losing the locality inherent to our approach, we introduced Transformers with a probabilistic loss function to make them more robust. However, in the framework of the geometric multigrid, we found that an enlarged patch size can be similarly effective.

We looked at the gap between the training objective and inference, noting that a lower validation loss does not necessarily translate into improved accuracy during a simulation. In fact, we found that larger corrections by the neural network can lead to instabilities, which manifest as localized, nonphysical features in the flow. To alleviate this, we proposed the application of replay buffers, a simple data augmentation technique that only requires the addition of a meta loop to retrain the neural networks on previous DNN-MG simulations. This allows one to improve the downstream accuracy and stability with great consistency, at the price of an increased training time and complexity.

We demonstrated DNN-MG's ability to generalize to different geometries and the doubling of the Reynolds number for a number of different metrics. These included time-local errors, which were decreased by a factor of up to 16 over the coarse solution, and mean velocities, which saw improvements of up to factor 5. Because the hybrid method can effectively leverage modern GPUs, the added runtime over the numerical method at the coarse resolution is negligible even with larger neural networks. In our experiments, we scaled up to \num{1e6} learnable weights, an order of magnitude more than there are degrees of freedom in our test cases.

% This could be of greater importance in transport dominated problems like the shallow-water equations, which we want to look at next.
In future work, we plan to look at transport dominated problems like the shallow-water equations. 
More work will also be needed to effectively utilize the corrected solution from previous steps as input to the neural network, which could be more important in this setting. Another interesting direction for future research would be to take a more probabilistic approach with diffusion models, especially for chaotic flows. Given the promising results of increased patch size, the effect of even larger patches could be explored in conjunction with neural network designs that exploit the internal structure of a patches to keep the computational cost manageable. Also, a more comprehensive study of the effects of the jump level would be valuable.

%%%%%%%%%%%%%%%%%%%%%%%%%%%%%%%%%%%%%%%%%%%%%%%%%%%%%%%%%%%%%%%%%%%%%%%%%%%%
\section*{Data availability}
The source code of DNN-MG, weights of the large neural networks and configurations to reproduce the experiments in this manuscript is available on Zenodo~\cite{Jendersie2026}.

\section*{CRediT authorship contribution statement}
\textbf{Robert Jendersie}: Conceptualization, Investigation, Methodology, Software, Visualization, Writing – original draft, Writing – review \& editing. \textbf{Nils Margenberg}: Methodology, Software, Writing – review \& editing. \textbf{Christian Lessig}:
Conceptualization, Funding Acquisition, Methodology, Supervision, Writing – original draft, Writing – review \& editing. \textbf{Thomas Richter}:
Conceptualization, Funding Acquisition, Methodology, Supervision, Writing – original draft, Writing – review \& editing.

\section*{Declaration of competing interest}
The authors declare that they have no known competing financial interests or personal relationships that could have appeared to influence the work reported in this paper.

\section*{Acknowledgment}
RJ and TR acknowledge support by Schmidt Sciences, Grant Number G-24-67790.
TR acknowledges the support of the GRK 2297 MathCoRe, funded by the Deutsche Forschungsgemeinschaft, Grant Number 314838170. TR and CL further acknowledge support of  the Deutsche Forschungsgemeinschaft, Grant Number 537063406.

%todo: any other project funding?

\appendix
%%%%%%%%%%%%%%%%%%%%%%%%%%%%%%%%%%%%%%%%%%%%%%%%%%%%%%%%%%%%%%%%%%%%%%%%%%%%
\section{Neural network details}\label{appendix:nn-details}
\begin{table}[!htb]
    \centering
%    \scriptsize
    \begin{tabular}{ccccccc}
        \toprule
        & small MLP & MLP & MLP-M1 & RNN & Transf \\
        \midrule
        hidden size $N_\text{hid}$ & 128 & 256 & 256 & 150 & 80 \\
        layers $L$ & 9 & 17 & 17 & 9 & 2 \\
        token size $N_\text{tok}$ & - & - & - & - & 128 \\
        ensemble size $N_\text{ens}$ & - & - & - & - & 4 \\ 
        encoder size & - & - & - & - & 256 \\
        encoder layers & - & - & - & - & 4 \\
        heads & - & - & - & - & 4 \\
        \midrule
        total weights & \num{126976} &\num{1011712} & \num{1044480} & \num{1045200} & \num{588672} \\
        \bottomrule
    \end{tabular}
    \caption{Size parameters for the neural networks used in the presented experiments. For Transf, $N_{\text{hid}}$ and $L$ describe each MLP in the tail ensemble.}
    \label{tab:nn-sizes}
\end{table}

%%%%%%%%%%%%%%%%%%%%%%%%%%%%%%%%%%%%%%%%%%%%%%%%%%%%%%%%%%%%%%%%%
\subsection{Multilayer perceptron}\label{appendix:nn-details:mlp}
An $\text{MLP}: \R^{N_\text{in}} \to \R^{N_\text{out}}$ is a composition of parameterized nonlinear functions
\begin{equation}\label{eq:mlp}
\text{MLP}(x) := (f_L\circ\cdots\circ f_1)(\vec{x})
\end{equation}
with layers $f_i : \R^{N_i} \to \R^{N_{i+1}}$, where
\begin{align*}
f_1(\vec{x}) &= \sigma (\text{LN}_1( W_1 \vec{x})) \\ %\in \mathcal{X}^1 = \mathcal{X}\\
f_j(\vec{x}) &= \sigma (\text{LN}_j(W_j \vec{x})) + \vec{x},\quad j = 2,\dots,\,L-1\, \\
f_L(\vec{x}) &= W_L \vec{x}.
\end{align*}
Each $f_i$ consists of a linear map with a learned weight matrix $W_i \in \R^{N_i \times N_{i+1}}$ followed by the layer normalization $N_i: \R^{n_{i+1}} \to \R^{n_{i+1}}$ which rescales each component of $\vec{x}$ according to
\begin{equation}\label{eq:layernorm}
\text{LN}_i(\vec{x}) = \frac{\vec{x} - E[\vec{x}]}{\sqrt{\text{Var}[\vec{x}] + 10^{-5}}} \odot \fvec{\gamma}_i + \fvec{\beta}_i.
\end{equation}
Here, $E[\vec{x}]$ and $\text{Var}[\vec{x}]$ are, respectively, the mean and variance of the components of the input vector $\vec{x}$, ``$\cdot \odot \cdot$'' is the Hadamard product between two tensors and $\fvec{\gamma}_i \in \R^{n_{i+1}}, \fvec{\beta}_i \in \R^{n_{i+1}}$ are the learnable parameters of an affine transformation. Afterward, the activation function $\sigma: \R \to \R$ is applied to each component. Finally, in the hidden layers, where the dimensions match, skip connections add the input $\vec{x}$ to the output.

%%%%%%%%%%%%%%%%%%%%%%%%%%%%%%%%%%%%%%%%%%%%%%%%%%%%%%%%%%%%%%%%%
\subsection{Recurrent neural network}\label{appendix:nn-details:rnn}
Our $\text{RNN}: \R^{N_\text{in}} \times \R^{(L-1) \times N_\text{hid}} \to \R^{N_\text{out}} \times \R^{(L-1) \times N_\text{hid}}$ is structurally similar to \cref{eq:mlp}. The linear maps are replaced by standard GRUs without bias for our RNN, with $\text{GRU}: \R^{N_\text{in}}\times\R^{N_\text{hid}} \to \R^{N_\text{hid}} \times\R^{N_\text{hid}}$. Consequently, the hidden state $\vec{h} \in \R^{(L-1)\times N_\text{hid}}$ needs to be tracked as an additional input and output for the $\text{RNN}: \R^{N_\text{in}} \times \R^{(L-1) \times N_\text{hid}} \to \R^{N_\text{out}} \times \R^{(L-1) \times N_\text{hid}}$ with
% (f_{L-1}^\text{out}\circ\cdots\circ f_1^\text{out})(x,h)
\begin{align}
    \text{RNN}^\text{out}(\vec{x},\vec{h}) &\coloneq  f_L(f_{L-1}^\text{out}(\dots f_1^\text{out}(\vec{x},\vec{h}_1), \vec{h}_{L-1})), \\
    \text{RNN}^\text{hid}(\vec{x},\vec{h}) & \coloneq \begin{bmatrix}f_1^{\text{hid}}(\vec{x},\vec{h}_1) & \dots & f_{L-1}^{\text{hid}}(\vec{x}_{L-1}, \vec{h}_{L-1})\end{bmatrix}
\end{align}
Each recurrent layer $f_i$ has its own hidden state that is zero-initialized but persists over the time-steps and produces an output $f_i^\text{out}$, as well as an updated hidden state $f_i^\text{hid}$, giving
\begin{align*}
f_1(\vec{x},\vec{h}) &= (\text{LN}_1( \text{GRU}^\text{out}(\vec{x},\vec{h})), \text{GRU}^\text{hid}(\vec{x},\vec{h})) \\
f_j(\vec{x},\vec{h}) &= (\text{LN}_j(\text{GRU}^\text{out}(\vec{x},\vec{h})) + \vec{x}, \text{GRU}^\text{hid}(\vec{x},\vec{h})),\quad j = 2,\dots,\,L-1\, \\
f_L(\vec{x}) &= W_L \vec{x}.
\end{align*}
For simplicity, layer normalization according to \cref{eq:layernorm} is applied only between GRUs and skip connections allow the inputs $\vec{x}$ to bypass the GRUs entirely. The final layer is just a linear projection that maps the latent vector to the output space.

%%%%%%%%%%%%%%%%%%%%%%%%%%%%%%%%%%%%%%%%%%%%%%%%%%%%%%%%%%%%%%%%%
\subsection{Transformer}\label{appendix:nn-details:transf}
Our transformer $\text{Transf}: \R^{N_{\text{seq}} \times N_\text{in}} \to \R^{N_\text{out}}$ consists of three sequential parts
\[
\text{Transf}(\vec{x}) := (\text{Tail} \circ \text{Enc} \circ \text{Emb}) (\vec{x}).
\]
First, the sequence of input patches $\vec{x} = \begin{bmatrix} \vec{x}_1 & \dots & \vec{x}_{N_s}\end{bmatrix}$ is processed by
\[
\text{Emb}(\vec{x}) = \begin{bmatrix} P \vec{x}_1 + \vec{p}_1 & \dots & P \vec{x}_{N_s} + \vec{p}_{N_s} & \vec{k}_\text{cls} \end{bmatrix}.
\]
Inputs are embedded, or projected to the internal token size $N_\text{tok}$, by the learned linear transform $P \in \R^{N_\text{tok} \times N_\text{in}}$ and a simple learned position encoding $p_i \in \R^{N_\text{tok}}$ is added to each token.  In order to force the neural network to collect relevant information from every token, another learned vector $\vec{k}_\text{cls} \in \R^{N_\text{tok}}$ is appended to the sequence which is later used to make the prediction.

The prepared sequence is fed into $\text{Enc}: \R^{N_{\text{seq}} \times {N_\text{tok}}} \to \R^{N_{\text{seq}}} \times \R^{N_\text{tok}}$, a standard BERT style transformer~\cite{devlin-etal-2019-bert} consisting of multiple encoder layers. We follow the common choice of GELU activations~\cite{hendrycks2023gaussianerrorlinearunits} and do not modify the overall structure.

Finally, the tail network $\text{Tail}: \R^{N_\text{seq} + 1} \times \R^{N_\text{tok}} \to \R^{N_\text{out}}$ selects the $\vec{k}_\text{cls}$ token and makes the prediction. Similar to \cite{Lessig2023}, the tail network is an ensemble of $N_\text{ens} \in \N $ shallow MLPs and the final result is computed as their mean, so that
\[
\text{Tail}(\vec{x}) = \frac{1}{N_\text{ens}} \sum_{j=1}^{N_\text{ens}} W_j^\text{out} \sigma (W_j^\text{hid} \vec{k}_\text{cls}).
\]
Each MLP has two linear layers defined by $W_j^\text{hid} \in \R^{N_\text{hid} \times N_\text{tok}}$ and $W_j^\text{out} \in \R^{N_\text{out} \times N_\text{hid}}$, as well as a single nonlinear GELU activation $\sigma$ in between.

%%%%%%%%%%%%%%%%%%%%%%%%%%%%%%%%%%%%%%%%%%%%%%%%%%%%%%%%%%%%%%%%%%%%%%%%%%%%
\section{Mesh details}\label{appendix:mesh-details}
\begin{table}[!htb]
    \centering
    \begin{tabular}{lcccp{4.7cm}p{5.5cm}}
    case & domain size & vertices & faces & obstacle position & obstacle radii \\
    \toprule
         ro1 & $(2.20, 0.41)$ & 16 &24 & $(0.20,0.20)$	& $(0.040, 0.060)$ \\
         ro3 & $(3.25, 0.41)$ & 28 & 44 & $(0.20,0.18)$ $(0.60, 0.22)$ $(1.00, 0.18)$ & $(0.045, 0.061)$ $(0.050, 0.050)$ $(0.040, 0.061)$ \\
         ro5 & $(3.25, 1.20)$ & 52 & 76 & $(0.20,0.22)$ $(1.00,0.18)$ $(0.60,0.60)$\newline$(0.20,0.98)$ $(1.00,1.02)$ & $(0.050, 0.050)$ $(0.050, 0.050)$ $(0.055, 0.035)$\newline$(0.050, 0.050)$ $(0.055, 0.035)$ \\
         ro6 & $(3.25, 1.2)$ & 56 & 83 & $(0.20, 0.18)$ $(0.60, 0.22)$ $(0.20, 0.58)$ \newline $(0.60, 0.62)$ $(0.20, 1.02)$ $(0.60, 0.98)$ & $(0.050, 0.050)$ $(0.050, 0.050)$ $(0.050, 0.050)$ \newline$(0.065, 0.055)$ $(0.050, 0.050)$ $(0.065, 0.055)$\\
         ro9 & $(3.25, 1.2)$ & 68 & 104 & $(0.20, 0.18)$ $(0.60, 0.22)$ $(1.00, 0.20)$ \newline $(0.20, 0.60)$ $(0.60, 0.58)$ $(1.00, 0.62)$ \newline $(0.20, 1.02)$ $(0.60, 0.98)$ $(1.00, 1.00)$ & $(0.049, 0.058)$ $(0.050, 0.050)$ $(0.050, 0.050)$ \newline$(0.050, 0.050)$ $(0.049, 0.058)$ $(0.050, 0.050)$ \newline$(0.049, 0.058)$ $(0.050, 0.050)$ $(0.050, 0.050)$\\
         ro3* & $(3.25, 0.41)$ & 28 & 44 & $(0.20, 0.18)$ $(0.60, 0.205)$ $(1.00, 0.18)$ & $(0.038, 0.055)$ $(0.050, 0.050)$ $(0.038, 0.055)$\\
         sq4 & $(2.45, 0.5)$ & 72 & 99 & $(0.30, 0.15)$ $(0.50, 0.15)$ $(0.30, 0.35)$ \newline $(0.50, 0.35)$ & $(0.050, 0.050)$ $(0.050, 0.050)$ $(0.050, 0.050)$ \newline$(0.050, 0.050)$\\
         sq6 & $(2.65, 0.7)$ & 112 & 149 & $(0.40, 0.15)$ $(0.60, 0.15)$ $(0.40, 0.35)$ \newline $(0.60, 0.35)$ $(0.40, 0.55)$ $(0.60, 0.55)$ & $(0.050, 0.050)$ $(0.050, 0.050)$ $(0.050, 0.050)$ \newline$(0.050, 0.050)$ $(0.050, 0.050)$ $(0.050, 0.050)$ \\
         sq9 & $(2.65, 0.7)$ & 112 & 158 & $(0.30, 0.15)$ $(0.50, 0.15)$ $(0.70, 0.15)$ \newline $(0.30, 0.35)$ $(0.50, 0.35)$ $(0.70, 0.35)$ \newline $(0.30, 0.55)$ $(0.50, 0.55)$ $(0.70, 0.55)$ & $(0.050, 0.050)$ $(0.050, 0.050)$ $(0.050, 0.050)$ \newline$(0.050, 0.050)$ $(0.050, 0.050)$ $(0.050, 0.050)$ \newline$(0.050, 0.050)$ $(0.050, 0.050)$ $(0.050, 0.050)$ \\
    \bottomrule
    \end{tabular}
    \caption{Details of the different base meshes used. The fine meshes are created by refining the round obstacle cases 5 times and the square obstacle cases 4 times. After refinment the obstacles are either elliptic (``ro$N$'') or square (``sq$N$'') holes in the mesh. In the modified ro3* the second obstacle is left as a square.}
    \label{tab:mesh-details}
\end{table}

\section{Experiment details}\label{appendix:experiment-details}
The large scale experiments used sq4, sq6, ro3, ro5 and ro9, simulated for 950 steps, as training sets. For the initial training data, we employed a noise schedule with strength $\sigma=0.02$ and probability determined by $c = 16$. The initial training also included an empty channel with a large cavity at the bottom. However, the simpler flow of the cavity case made it less interesting and since its impact was negligible it was not included for replay generation. In total, removing samples from the 50 warmup steps, this yielded \num{3.6e7} unique samples for initial training. After truncating the simulations to 650 steps for the replays, roughly \num{4.5e7} samples where gathered for one training epoch. As validation set we used sq9, ro6, ro1 and ro1 with noise augmentation ($\sigma=0.04, c=1$) to select the best neural network during training, resulting in \num{1.2e7} samples. Additional sets with other meshes where evaluated during validation and thus counted in the training time, but did not factor into the selection.
%36443648, 44925952

\begin{figure}[!htb]
    \centering
    \ifbuild
    \tikzsetnextfilename{velocity_artifacts}
    \begin{tikzpicture}
        \node(ro9){\includegraphics[width=0.95\textwidth]{imgs/gc9_0_ifix_mlp_large_62_step450.png}};
        \node[below=0cm of ro9]{\includegraphics[width=0.795\textwidth]{imgs/c9_5e-3r_0_ifix_mlp_large_62_step197_gc9scale.png}};
    \end{tikzpicture}
    \else
    \includegraphics{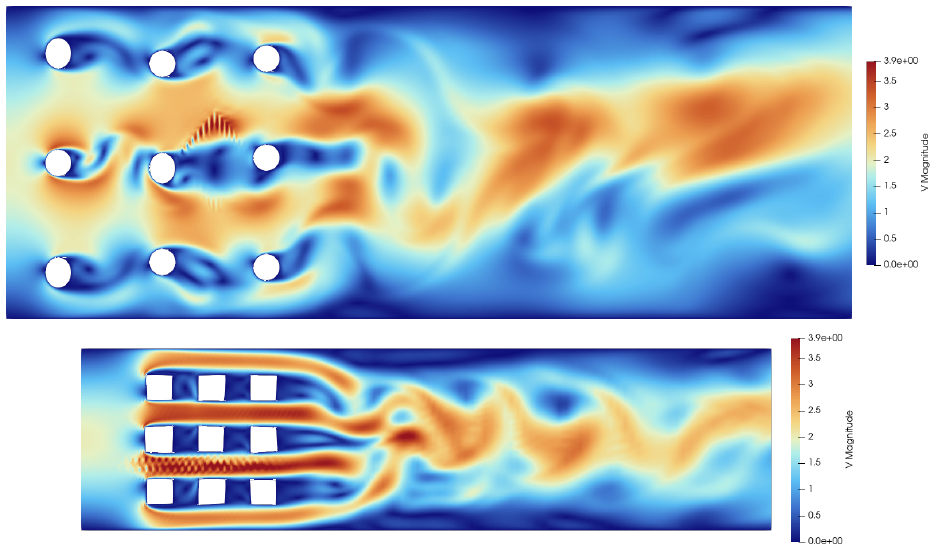}
    \fi
    \caption{Velocity magnitude snapshots from two simulations with DNN-MG that clearly show artifacts. In both cases large MLP without replay training was used. In ro9 (top) at $t=4.5$, high frequency waves are visible behind the middle obstacle that closely follow the mesh geometry. In sq9* at $t=1.97$, significant perturbations develop in the lower channel between obstacles. Furthermore, there are higher-frequency wave fronts that travel diagonal to the channel behind the obstacles. These fronts are more visible in motion.}
    \label{fig:artifacts-velocity}
\end{figure}

\begin{figure}[!htb]
    \centering
    \ifbuild
\else
    \includegraphics{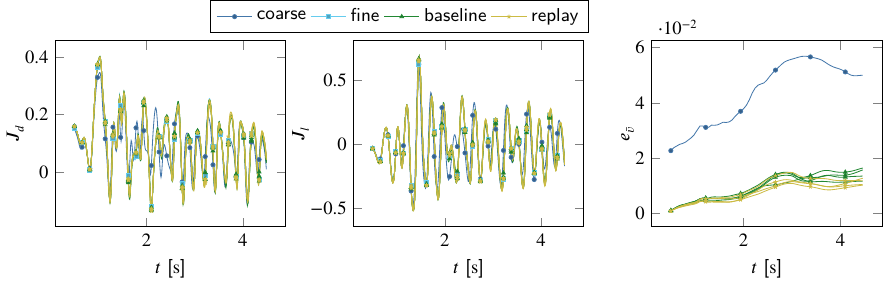}
    \fi
	\caption{Metrics over time for simulations of ro9 which are summarized in \cref{tab:replay-gc9}. The mean velocity error $e_{\bar{v}}$~\eqref{eq:mean-vel-error} is taken over a growing time-interval, keeping the start time fixed and increasing the number of time-steps.}
	\label{fig:replay-gc9}
\end{figure}

\begin{figure}
    \centering
    \ifbuild
\else
	\includegraphics{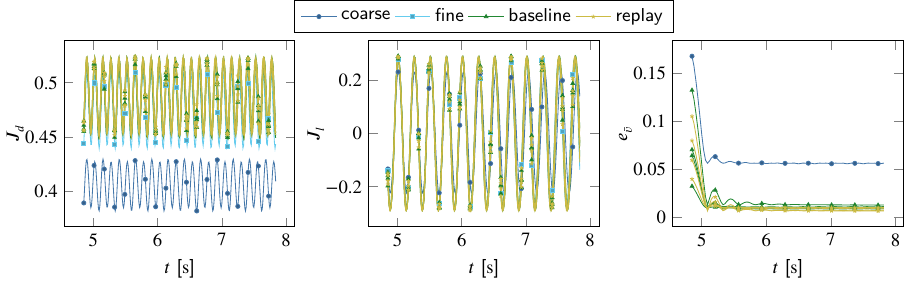}
	\fi
	\caption{Metrics over time for simulations of ro1 which are summarized in \cref{tab:replay-gc1}. The mean velocity error $e_{\bar{v}}$~\eqref{eq:mean-vel-error} is taken over a growing time-interval, keeping the start time fixed and increasing the number of time-steps. }
	\label{fig:replay-gc1}
\end{figure}

%%%%%%%%%%%%%%%%%%%%%%%%%%%%%%%%%%%%%%%%%%%%%%%%%%%%%%%%%%%%%%%%%
\begin{figure}
    \centering
    \ifbuild
\else
	\includegraphics{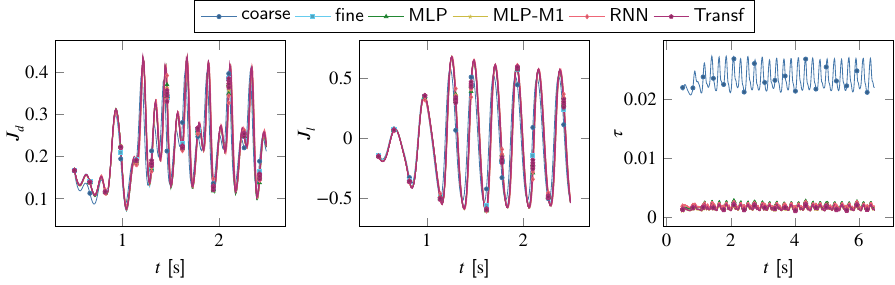}
	\fi
	\caption{Metrics over time for simulations of ro6 which are summarized in \cref{tab:architectures-gc6}.}
	\label{fig:architectures-gc6}
\end{figure}

\begin{figure}
    \centering
    \ifbuild
\else
	\includegraphics{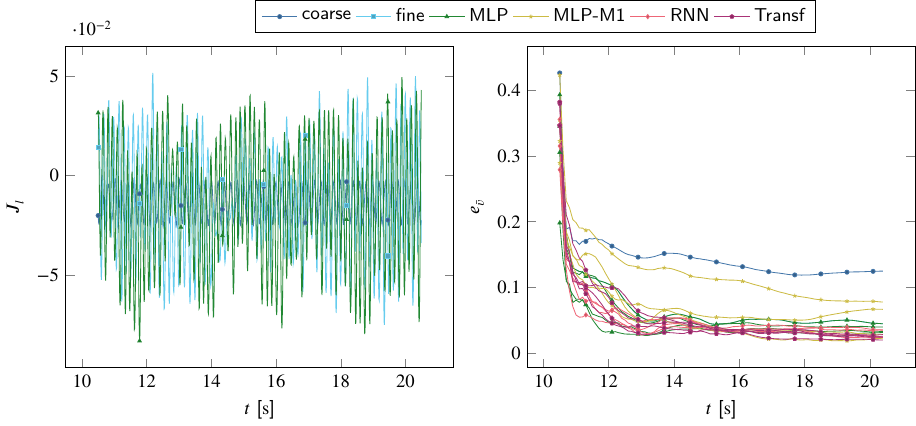}
	\fi
	\caption{Metrics over time for simulations of sq9* which are summarized in \cref{tab:architectures-c9}. The lift functional $J_l$~\eqref{eq:lift} is shown only for a single MLP for clarity but the other neural networks follow a similar trend. The mean velocity error $e_{\bar{v}}$~\eqref{eq:mean-vel-error} is taken over a growing time-interval, keeping the start time fixed and increasing the number of time-steps.}
	\label{fig:architectures-c9}
\end{figure}

\begin{figure}
    \centering
    \ifbuild
\else
	\includegraphics{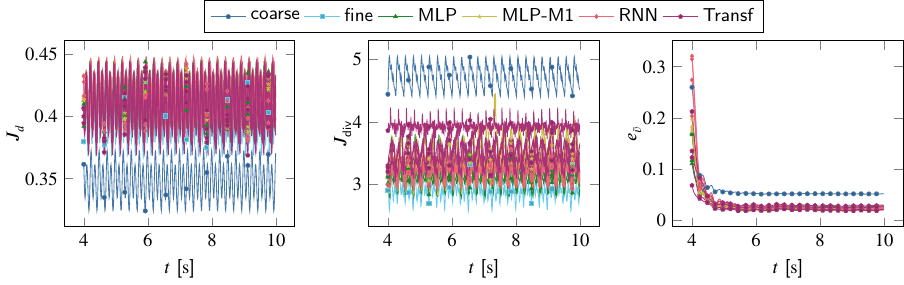}
	\fi
	\caption{Metrics over time for simulations of ro3* which are summarized in \cref{tab:architectures-sgc3}.}
	\label{fig:architectures-sgc3}
\end{figure}

%%%%%%%%%%%%%%%%%%%%%%%%%%%%%%%%%%%%%%%%%%%%%%%%%%%%%%%%%%%%%%%%%%%%%%%%%%%%
\section{Fine velocity inputs}\label{appendix:fine-vel-inputs}
To use the NN-corrected velocity field in future time-steps, modifications to both the training and the solver are necessary. In the basic implementation, we construct training samples with the reference solution as velocity inputs and only mask it with the interpolated coarse solution for the current time-step. In the solver, the current neural network prediction is added right back to the input buffer. Additional care is needed in the solver when inputs or outputs are normalized, since the velocities have to be scaled and shifted appropriately before adding them back.

For transformers trained with this setup, we observed reductions in the final validation loss by a factor of 2 on round obstacle cases like ro1 and ro6. However, during the simulation the DNN-MG solution would closely follow the coarse solution, with the neural network only making tiny corrections. This is because a neural network trained this way has not learned how to reconstruct the higher fidelity field. To alleviate this, we modified the training to use input velocities that are randomly interpolated between the coarse and the fine solution as
\[
\vec{v}_{\text{inp}} = t\intsol{\vec{v}} + (1-t)\refsol{\vec{v}},\quad t = \min(1, 4r),
\]
where $r \in [0,1]$ is sampled from a uniform distribution. 
%Nevertheless, even weighted towards the fine solution, the randomized training yielded significantly less improvement in the final loss in our experiments compared to just the fine solution. 
%However the convergence was still faster than with just the coarse velocity. 
During inference, the neural networks trained on coarse-fine interpolated inputs performed better initially than the fine velocity ones, but still worse than the baseline. Furthermore, the direct pathway for feedback caused significant stability issues that would often cause the numerical simulation to explode after a few hundred steps.

To improve the stability of fine velocity inputs, we tried replay training as discussed in \cref{sec:replay}. As this did not proof sufficient, we tried finetuning with back-propagation through time. Adding the predictions back to the inputs during training opens a pathway for gradient computations over multiple steps that directly addresses the feedback loop albeit the contribution from the solver is missing. We tried a roll-out for 8 steps that made the training roughly 5 times slower. While this helped, the neural network still performed worse in the hybrid simulations than neural networks trained only on coarse inputs.
%However this did not prove effective and in fact made the neural network worse. 
One problem that remains is the use of mini-batches constructed from random patches which disregards the interaction between neighboring patches. We tried to solve this by always training with pairs of neighboring patches but the training signal from just one neighbor is likely to weak. Another problem in our BPTT training is that we used precomputed multi-step training data and only updated the velocity, causing a mismatch in the residuals and a greater distributional shift between training and inference. A proper roll-out involving the solver itself on a full domain and complete gradients would probably be more effective.

%%%%%%%%%%%%%%%%%%%%%%%%%%%%%%%%%%%%%%%%%%%%%%%%%%%%%%%%%%%%%%%%%%%%%%%%%%%%
\section{Jump level}\label{appendix:jmplvl}
{As supplement to \cref{sec:generalization} we provide data from 4 large MLPs each with an increased jump level. To make them comparable, the coarse level is reduced by one, henceforth denoted as coarse $(-2)$, such that the NNs are tasked to reconstruct the same solution but with less information from the numerical solver. We present 2 variants, MLP-J2 ($J=2, M=0$) which uses patches corresponding to cells on the fine mesh and MLP-MJ2M1 ($J=2, M=1$) which has patches corresponding to the former coarse $(-1)$ mesh. This task is significantly more difficult as the solutions on the coarse $(-2)$ meshes have only 576 (ro1) to 2667 (ro6) nodes which is insufficient for the typical dynamics to develop.

The new coarse solution $(-2)$ is significantly worse, and consequently the NNs perform strictly worse than the MLPs at $J=1$. In particular, coarse $(-2)$ significantly underestimates the forces acting on the round obstacles and the NNs tend to overcompensate lift (\cref{tab:jmplvl-gc6}, \cref{fig:jmplvl-gc6}) or drag (\cref{fig:jmplvl-gc1}). The local error and divergence are still very close to the fine reference for all cases and constitute a significant improvement over the numerical solution coarse $(-1)$~(\cref{tab:jmplvl-gc6}, \cref{tab:jmplvl-c9}, \cref{tab:jmplvl-sgc3}). The mean velocity error over longer time periods is around coarse $(-1)$. On ro9*, the metric is of limited value here as we barely see any convergence from coarse $(-2)$ to coarse $(-1)$~(\cref{tab:jmplvl-c9}). On ro3*, the NNs manage to bridge most of the gap to the fine reference and MLP-J2M1 outperforms the numerical solver coarse $(-1)$~(\cref{tab:jmplvl-sgc3}).

%%%%%%%%%%%%%%%%%%%%%%%%%%%%%%%%%%%%%%%%%%%%%%%%%%%%%%%%%%%%%%%%%
\begin{table}
    \centering
    \begin{tabular}{lcccccccc}
    \toprule
     & \multicolumn{2}{c}{mean $e_{J_d}$ \eqref{eq:drag-error}} & \multicolumn{2}{c}{mean $e_{J_l}$ \eqref{eq:lift-error}} & \multicolumn{2}{c}{mean $J_{\text{div}}$ \eqref{eq:div}} & \multicolumn{2}{c}{mean $\tau$ \eqref{eq:trunc-error}} \\
     & mean & best & mean & best & mean & best & mean & best \\
    \midrule
    coarse (-2) & 0.077 & 0.077 & 0.310 & 0.310 & 6.275 & 6.275 & 0.366 & 0.366 \\
    coarse (-1) & \textit{0.031} & \textit{0.031} & \textit{0.098} & \textit{0.098} & 3.307 & 3.307 & 0.024 & 0.024 \\
    MLP & \textbf{0.008} & \textbf{0.006} & \textbf{0.025} & \textbf{0.019} & \textbf{1.348} & \textbf{1.337} & \textbf{0.002} & \textbf{0.002} \\
    MLP-J2 & 0.045 & 0.038 & 0.126 & 0.115 & 2.124 & 1.993 & 0.008 & 0.008 \\
    MLP-J2M1 & 0.050 & 0.048 & 0.147 & 0.125 & \textit{1.634} & \textit{1.630} & \textit{0.005} & \textit{0.005} \\
    \midrule
    fine & 0.000 & 0.000 & 0.000 & 0.000 & 1.253 & 1.253 & 0.000 & 0.000 \\
    \bottomrule
    \end{tabular}
    \caption{Simulation results on ro6 for 4 NNs each for different jump levels. All metrics are computed starting right after the warmup at $n_0=50$. Drag error $e_{J_d}$ and lift error $e_{J_l}$ are taken over $k=200$ steps to keep the influence of the temporal shift small. Divergence $J_{\text{div}}$ and the local error $\tau$ are averaged over $k=600$ steps.. The entries MLP, coarse (-1) and fine are the same as \cref{tab:architectures-gc6}. The values are plotted over time in \cref{fig:jmplvl-gc6}.}
	\label{tab:jmplvl-gc6}
\end{table}
\begin{figure}
    \centering
    \ifbuild
\else
	\includegraphics{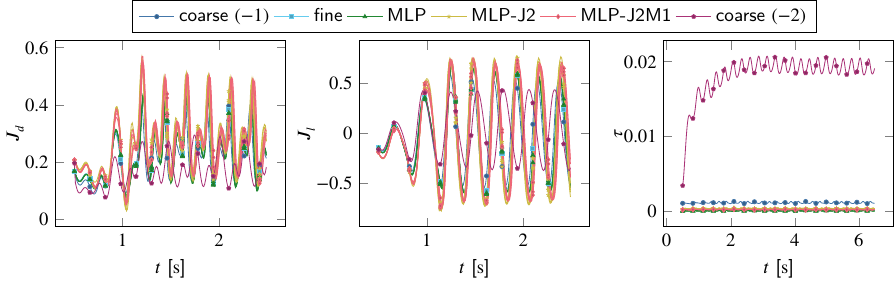}
	\fi
	\caption{Metrics over time for simulations of ro6 which are summarized in \cref{tab:jmplvl-gc6}.}
	\label{fig:jmplvl-gc6}
\end{figure}

%%%%%%%%%%%%%%%%%%%%%%%%%%%%%%%%%%%%%%%%%%%%%%%%%%%%%%%%%%%%%%%%%
\begin{table}
    \centering
    \begin{tabular}{lcccccccc}
    \toprule
     & \multicolumn{2}{c}{$\min J_l$ \eqref{eq:lift}} & \multicolumn{2}{c}{$\max J_l$ \eqref{eq:lift}} & \multicolumn{2}{c}{mean $J_{\text{div}}$ \eqref{eq:div}} & \multicolumn{2}{c}{$e_{\bar{v}}$ \eqref{eq:mean-vel-error}} \\
     & mean & best & mean & best & mean & best & mean & best \\
    \midrule
    coarse (-2) & -0.140 & -0.140 & 0.116 & 0.116 & 9.969 & 9.969 & 0.107 & 0.107 \\
    coarse (-1) & -0.026 & -0.026 & \textit{-0.000} & -0.000 & 7.319 & 7.319 & 0.094 & 0.094 \\
    MLP & \textbf{-0.079} & \textbf{-0.075} & \textbf{0.051} & \textbf{0.052} & \textbf{5.122} & \textbf{5.108} & \textbf{0.029} & \textbf{0.025} \\
    MLP-J2 & \textit{-0.047} & \textit{-0.056} & -0.008 & \textit{0.007} & 5.437 & 5.420 & \textit{0.093} & \textit{0.078} \\
    MLP-J2M1 & -0.027 & -0.030 & -0.011 & -0.009 & \textit{5.248} & \textit{5.243} & 0.111 & 0.093 \\
    \midrule
    fine & -0.075 & -0.075 & 0.051 & 0.051 & 5.125 & 5.125 & 0.000 & 0.000 \\
    \bottomrule
    \end{tabular}
    \caption{Simulation results on sq9* for 4 NNs each for different jump levels. All metrics are taken from $n_0=1050$ for $k=1000$ steps. The entries MLP, coarse (-1) and fine are the same as \cref{tab:architectures-c9}. The values are plotted over time in \cref{fig:jmplvl-c9}.}
	\label{tab:jmplvl-c9}
\end{table}
\begin{figure}
    \centering
    \ifbuild
\else
	\includegraphics{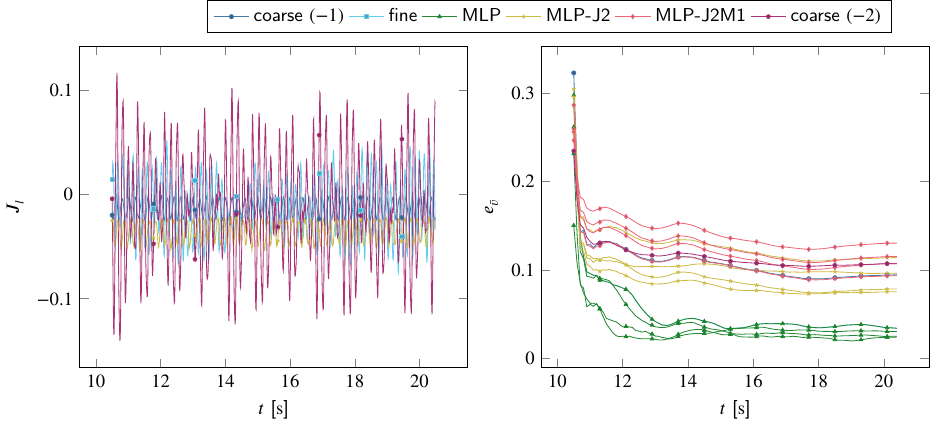}
	\fi
	\caption{Metrics over time for simulations of sq9* which are summarized in \cref{tab:jmplvl-c9}. The lift functional $J_l$~\eqref{eq:lift} is shown only for a single MLP-J2 for clarity but the other neural networks follow a similar trend. The mean velocity error $e_{\bar{v}}$~\eqref{eq:mean-vel-error} is taken over a growing time-interval, keeping the start time fixed and increasing the number of time-steps.}
	\label{fig:jmplvl-c9}
\end{figure}

%%%%%%%%%%%%%%%%%%%%%%%%%%%%%%%%%%%%%%%%%%%%%%%%%%%%%%%%%%%%%%%%%
\begin{figure}
\centering
\ifbuild
\else
\includegraphics{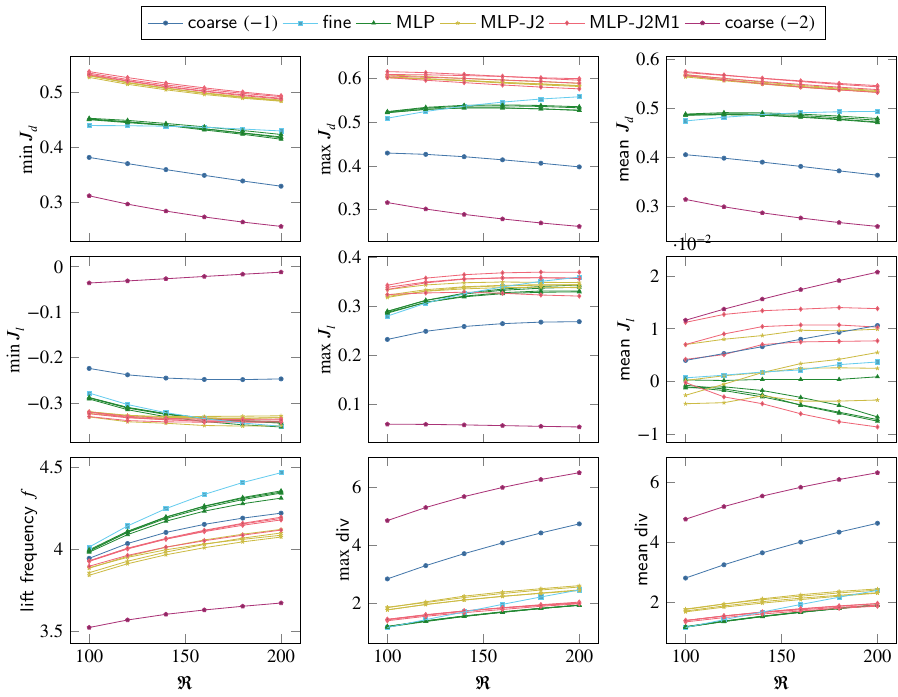}
\fi
	\caption{Functionals on ro1 for 4 NNs each for different jump levels, taken from $n_0=850$ over $k=200$ steps for simulations with increasing Reynolds number $\Rey$. The viscosity $\nu = 0.001$ is the same in all training cases, giving $\Rey=100$, and reduced to $\nu=0.0005$ for the highest Reynolds number.}
	\label{fig:jmplvl-gc1}
\end{figure}

%%%%%%%%%%%%%%%%%%%%%%%%%%%%%%%%%%%%%%%%%%%%%%%%%%%%%%%%%%%%%%%%%
\begin{table}
    \centering
    \begin{tabular}{lcccccccc}
    \toprule
     & \multicolumn{2}{c}{mean $J_d$ \eqref{eq:drag}} & \multicolumn{2}{c}{amp $J_d$ \eqref{eq:drag}} & \multicolumn{2}{c}{mean $J_{\text{div}}$ \eqref{eq:div}} & \multicolumn{2}{c}{$e_{\bar{v}}$ \eqref{eq:mean-vel-error}} \\
     & mean & best & mean & best & mean & best & mean & best \\
    \midrule
    coarse (-2) & 0.290 & 0.290 & 0.013 & 0.013 & 8.277 & 8.277 & 0.206 & 0.206 \\
    coarse (-1) & \textit{0.348} & \textit{0.348} & \textit{0.049} & \textit{0.049} & 4.726 & 4.726 & 0.041 & 0.041 \\
    MLP & \textbf{0.413} & \textbf{0.412} & \textbf{0.066} & \textbf{0.065} & \textbf{3.188} & \textbf{3.097} & \textbf{0.020} & \textbf{0.017} \\
    MLP-J2 & 0.517 & 0.513 & 0.084 & 0.079 & 4.299 & 4.028 & 0.058 & 0.052 \\
    MLP-J2M1 & 0.512 & 0.512 & 0.078 & 0.077 & \textit{4.248} & \textit{3.962} & \textit{0.033} & \textit{0.030} \\
    \midrule
    fine & 0.399 & 0.399 & 0.062 & 0.062 & 2.919 & 2.919 & 0.000 & 0.000 \\
    \bottomrule
    \end{tabular}
\caption{Simulation results on ro3* for 4 NNs each for different jump levels. All metrics are computed starting at $n_0=400$ and computed over $k=600$ steps. The amplitude of the drag is computed as $\text{amp} J_d = \max J_d - \min J_d$. The entries MLP, coarse (-1) and fine are the same as \cref{tab:architectures-sgc3}. The values are plotted over time in \cref{fig:jmplvl-sgc3}.}
	\label{tab:jmplvl-sgc3}
\end{table}
\begin{figure}
    \centering
    \ifbuild
\else
	\includegraphics{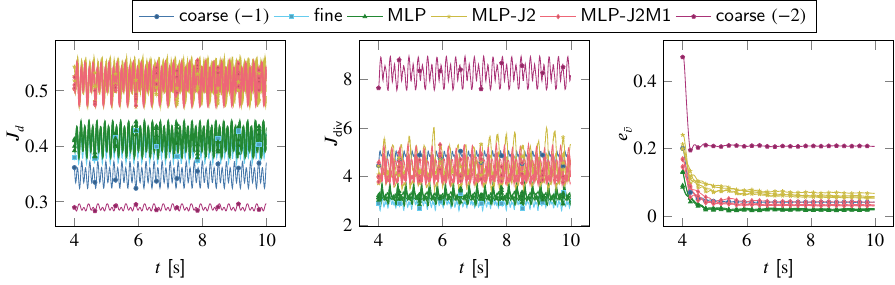}
	\fi
	\caption{Metrics over time for simulations of ro3* which are summarized in \cref{tab:jmplvl-sgc3}.}
	\label{fig:jmplvl-sgc3}
\end{figure}

%%%%%%%%%%%%%%%%%%%%%%%%%%%%%%%%%%%%%%%%%%%%%%%%%%%%%%%%%%%%%%%%%%%%%%%%%%%%
%%%%%%%%%%%%%%%%%%%%%%%%%%%%%%%%%%%%%%%%%%%%%%%%%%%%%%%%%%%%%%%%%%%%%%%%%%%%
% bibliography
\bibliographystyle{elsarticle-num} 
\bibliography{climate.bib,neuralnets.bib,paper_biber.bib,lit.bib}

\end{document}